# Tidy subgroups for commuting automorphisms of totally disconnected groups: an analogue of simultaneous triangularisation of matrices

George A. Willis

ABSTRACT. Let $\alpha$ be an automorphism of the totally disconnected group $G$. The compact open subgroup, $V$, of $G$ is *tidy* for $\alpha$ if $[\alpha(V') : \alpha(V') \cap V']$ is minimised at $V$, where $V'$ ranges over all compact open subgroups of $G$. Identifying a subgroup tidy for $\alpha$ is analogous to identifying a basis which puts a linear transformation into Jordan canonical form. This analogy is developed here by showing that commuting automorphisms have a common tidy subgroup of $G$ and, conversely, that a group $\mathfrak{H}$ of automorphisms having a common tidy subgroup $V$ is abelian modulo the automorphisms which leave $V$ invariant. Certain subgroups of $G$ are the analogues of eigenspaces and corresponding real characters of $\mathfrak{H}$ the analogues of eigenvalues.

## 1. Introduction

Let $G$ be a totally disconnected locally compact group and let $\alpha$ be a continuous automorphism of $G$. Then $G$ has a base of neighbourhoods of the identity consisting of compact open subgroups, see [4, Theorem 7.7], and for each such subgroup $V$ the index $[\alpha(V) : \alpha(V) \cap V]$ is finite. The *scale* of $\alpha$ is the positive integer

(1) $\quad s(\alpha) = \min\{[\alpha(V) : \alpha(V) \cap V] : V \text{ is a compact open subgroup of } G\}$

and those compact open subgroups $V$ where this minimum is attained are said to be *tidy* for $\alpha$. Tidy subgroups and the scale function have played an essential role in the proofs of several conjectures about totally disconnected groups. For instance, it is shown in [6] that there are no totally disconnected locally compact groups satisfying the 'strangeness' condition defined [5] and, in [10], the scale function alone is used to show that the set of periodic elements in a totally disconnected group is closed (which is not true of connected groups).

In the special case when $G$ is the additive group of an $n$-dimensional $p$-adic linear space, choosing a basis identifies $G$ with $\mathbf{Q}_p^n$ and identifies automorphisms of $G$ with elements of $GL(n, \mathbf{Q}_p)$. Suppose that the automorphism $\alpha$ is in Jordan canonical form and let its diagonal entries be $a_j$, $j = 1, \ldots, n$, with $p$-adic absolute

Research supported by A.R.C. Grant A69700321.

values $|a_j| = p^{-k_j}$. Then the compact open subgroup $\mathbf{Z}_p^n$ is tidy for $\alpha$ and $s(\alpha) = \prod_{k_j \leq 0} p^{-k_j}$. Tidy subgroups and the scale function of $\alpha$ are thus described in terms of a basis which makes $\alpha$ triangular. This description may be shown directly in this example but it can also be shown with the aid of the following general characterisation of tidy subgroups.

The compact open subgroup $V$ of $G$ is tidy for $\alpha$ if and only if it satisfies

**T1**($\alpha$): $V = V_+ V_-$, where $V_\pm := \bigcap_{n \geq 0} \alpha^{\pm n}(V)$; and
**T2**($\alpha$): $V_{++} := \bigcup_{n \geq 0} \alpha^n(V_+)$ is closed.

In fact, the original definition of tidy subgroup, made in [9], was that $V$ should satisfy **T1**($\alpha$) and **T2**($\alpha$) and [8, Theorem 3.1] shows that this definition is equivalent to the one given above. It is immediate from their definition that $V_+$ and $V_-$ are compact subgroups of $V$, that $\alpha(V_+) \geq V_+$ and that $\alpha(V_-) \leq V_-$. Condition **T1**($\alpha$) thus asserts that tidy subgroups are the product of a subgroup on which $\alpha$ expands and a subgroup on which $\alpha$ contracts and this provides a local description of $\alpha$. The index $[\alpha(V_+) : V_+]$ is the expansion factor on $V_+$ and is equal to the scale of $\alpha$. Condition **T2**($\alpha$) is equivalent to the more global assertion that $\alpha$-orbits cannot leave and re-enter $V$, [9, Lemma 3].

Return now to the example of the automorphism $\alpha$ of $\mathbf{Q}_p^n$ having the tidy subgroup $V = \mathbf{Z}_p^n$. Direct calculation shows that

$$V_+ = \bigoplus_{k_j \leq 0} \mathbf{Z}_p, \qquad V_- = \bigoplus_{k_j \geq 0} \mathbf{Z}_p \text{ and}$$

$$V_{++} = \bigoplus_{k_j < 0} \mathbf{Q}_p \oplus \bigoplus_{k_j = 0} \mathbf{Z}_p,$$

so that these subgroups too may be described in terms of any basis for which $\alpha$ is in Jordan canonical form. The notion of subgroup tidy for an automorphism is in this way analogous to the Jordan canonical form of a matrix. In cases when $\alpha$ cannot be triangularised without extending the base field, tidy subgroups can be described in terms of factors of the characteristic polynomial of $\alpha$.

Tidy subgroups are seen to be a sort of canonical form in several other classes of groups as well, including: Lie groups over local fields, [2, 3]; automorphism groups of homogeneous trees; [9, §3]; and certain restricted products of finite groups, [7].

A set, $\mathfrak{A}$, of automorphisms of $G$ will be said to have a *common tidy subgroup* if there is a compact open subgroup $V$ of $G$ which is tidy for every $\alpha$ in $\mathfrak{A}$. The example when $G = \mathbf{Q}_p^n$ and automorphisms are linear transformations suggests that any commuting set of automorphisms should have a common tidy subgroup. This is shown for finite sets of automorphisms in §3 and for finitely generated abelian groups of automorphisms in §5. In the converse direction, it shown in §4 that, if $\mathfrak{H}$ is a group of automorphisms of $G$ having a common tidy subgroup $V$, then $\mathfrak{H}$ is abelian modulo the subgroup of automorphisms which leave $V$ invariant. (There are pairs of automorphisms, $\alpha$ and $\beta$, having a common tidy subgroup such that $\langle \alpha, \beta \rangle$ does not have a common tidy subgroup and so a converse applying to sets of automorphisms rather than groups cannot be expected.)

An important part of the proofs of the results in §4 and §5 is a refinement of the property **T1** of tidy subgroups. It is shown that, if $V$ is tidy for automorphisms $\alpha_1, \ldots, \alpha_n$, then $V$ is a product of subgroups, $V_\mathfrak{a}$, such that each $\alpha_j$ is expanding, contracting or neither on each $V_{\mathfrak{a}}$. In view of the example when $G = \mathbf{Q}_p^n$, each of the subgroups $V_{\mathfrak{a}}$ may be regarded as analogous to a common eigenspace for the automorphisms $\alpha_j$ and the expansion or contraction factor of each $\alpha_j$ as (an absolute value of) an eigenvalue. Properties of these 'eigenspaces' and 'eigenvalues' are worked out in §6.

The arguments used here rely on facts about tidy subgroups and the scale established in earlier papers. In particular, [9, Corollary 3] shows that the scale function $s : \mathrm{Aut}(G) \to \mathbf{Z}^+$ satisfies:

**S1:** $s(\alpha) = 1 = s(\alpha^{-1})$ if and only if there is a compact open subgroup $U$ of $G$ with $\alpha(U) = U$;

**S2:** $s(\alpha^n) = s(\alpha)^n$ for every $n \geq 0$; and

**S3:** $\Delta(\alpha) = s(\alpha)/s(\alpha^{-1})^{-1}$, where $\Delta : \mathrm{Aut}(G) \to \mathbf{Q}^+$ is the modular function.

Also, a procedure for producing tidy subgroups is given in [9] and improved in [8]. Yet another variant of this procedure is used to produce subgroups tidy for several commuting automorphisms in §3. The following section describes this variation and shows that it produces the same tidy subgroup as the version in [8].

## 2. The Tidying Procedure

Tidy subgroups are produced by choosing an arbitrary compact open subgroup $U$ and then carrying out three steps. Here is the procedure for tidying $U$ as described in [8].

**Step 1:** *'trim' $U$ to achieve property $\mathbf{T1}(\alpha)$.* [9, Lemma 1] shows that there is an $n$ such that $\bigcap_{k=0}^{n} \alpha^k(U)$ satisfies $\mathbf{T1}(\alpha)$. Choose such an $n$ and put

$$V := \bigcap_{k=0}^{n} \alpha^k(U).$$

**Step 2:** *locate the obstruction to $V_{++}$ being closed.* Define

(2) $$\mathcal{L} := \{l \in G : \alpha^j(l) \notin V \text{ for only finitely many } n\}$$

and put $L := \mathcal{L}^-$. Then $L$ is compact by [9, Lemma 6] and $V$ satisfies $\mathbf{T2}(\alpha)$ if and only if $L \leq V$.

**Step 3:** *join $V$ and $L$ to obtain a tidy subgroup.* Define

(3) $$V' := \{v \in V : lvl^{-1} \in VL \text{ for all } l \in L\}.$$

Then $W := V'L$ is a compact open subgroup of $G$ and satisfies the tidiness criteria $\mathbf{T1}(\alpha)$ and $\mathbf{T2}(\alpha)$, see [8, Lemmas 3.2–3.10].

The point in Step 3 is that $VL$ need not be a group but $V'L$ is. It is important here, and shown in Lemma 2.4, that $V'L$ is in fact maximal among groups contained in the set $VL$. Step 3 is the only place where the tidying procedures in [9] and [8] differ: that in [9] may produce a smaller tidy subgroup.

The new procedure differs from the previous ones at Step 2. Suppose that we have defined $V$ as in Step 1.

**Step 2a:** *locate a smaller obstruction to $V_{++}$ being closed.* Define

(4) $$\mathcal{K} := \{l \in G : \lim_{j \to \infty} \alpha^j(l) = e \text{ and } \{\alpha^{-j}(l)\}_{j \geq 0} \text{ is bounded}\}$$

and put $K := \mathcal{K}^-$.

**Step 3a:** *join V and K*. Define

(5) $$V'' := \{v \in V : lvl^{-1} \in VK \text{ for all } l \in K\}.$$

Now we will show that $V''K$ is a compact open subgroup of $G$ which is tidy for $\alpha$ by showing that it is in fact the same group $W$ produced by the earlier procedure. The advantage gained by using $K$ rather than $L$ derives from the fact that the definition of $K$ is independent of the subgroup $U$ chosen at the beginning.

In the remainder of this section and the next $V$, $V'$, $V''$, $L$, $K$ and $W$ will be the groups defined in the above steps.

**Proposition 2.1.** *The set $V''K$ is a compact open subgroup of $G$ and is the same tidy subgroup $W$ as constructed in the previous procedure.*

The proof of the proposition rests on some lemmas.

**Lemma 2.2.** *$K$ is a compact normal subgroup of $L$.*

**Proof.** It is easily checked that $\mathcal{K}$, and hence $K$, is a group. To show that $K$ is contained in $L$ it suffices to show that for each $k$ in $\mathcal{K}$ and each compact open subgroup $B$ of $G$ there is an element, $l$, of $\mathcal{L}$ in $kB$. Since $V$ is open it may be supposed that $B \leq V$ and since, by [9, Lemma 1], each compact open subgroup contains one satisfying $\mathbf{T1}(\alpha)$, it may be supposed that $B$ satisfies $\mathbf{T1}(\alpha)$.

Let $k$ be in $\mathcal{K}$ and $B$ be such a compact open subgroup. The definition of $\mathcal{K}$ implies that we may choose

(6)     a positive integer $m$ such that $\alpha^n(k)$ is in $B$ for all $n \geq m$

and

(7)     positive integers $i$ and $j$, with $i - j > m$, such that $\alpha^{-i}(k) \in B\alpha^{-j}(k)$.

It follows from (7) that $k\alpha^{i-j}(k^{-1})$ is in $\alpha^i(B)$. Then, using that $B$ satisfies $\mathbf{T1}(\alpha)$, $k\alpha^{i-j}(k^{-1})$ may be factored as

(8)     $k\alpha^{i-j}(k^{-1}) = lu$, where $\alpha^n(l) \in B$ for $n \leq -i$ and $\alpha^n(u) \in B$ for $n \geq -i$.

Since $i - j$ is positive, (6) implies that $\alpha^n(k\alpha^{i-j}(k^{-1}))$ is in $B$ whenever $n \geq m$. Together with (8), this implies that

$$\alpha^n(l) \in B \text{ for } n \leq -i \text{ or } n \geq m.$$

Therefore, since $B$ is contained in $V$, $l$ belongs to $\mathcal{L}$. To see that $l$ belongs to $kB$, note that the element $u$ found in (8) belongs to $B$ and that, since $i - j > m$, (6) gives that $k\alpha^{i-j}(k^{-1})$ is in $kB$.

Lemma 6 in [9] shows that $L$ is compact. Hence $K$ is compact too.

Now let $l$ be in $L$ and $k$ in $\mathcal{K}$. Then $\alpha^j(l)$ is in $L$ for all $j$ and so

$$\lim_{j \to \infty} \alpha^j(lkl^{-1}) = \lim_{j \to \infty} \alpha^j(l)\alpha^j(k)\alpha^j(l^{-1}) = e$$

because $L$ is compact. It follows similarly that $\{\alpha^{-j}(lkl^{-1})\}_{j \geq 0}$ is bounded and so $lkl^{-1}$ belongs to $\mathcal{K}$. Hence $\mathcal{K}$ is normal in $L$ and so $K$ is normal too. □

**Lemma 2.3.** *$L = (V \cap L)K$.*

**Proof.** The inclusion $(V \cap L)\mathcal{K} \subseteq L$ is obvious. Hence it suffices to find an element of $\mathcal{K}$ in each $(V \cap L)$-coset of $L$.

Suppose for now that $L$ is separable. Choose a sequence, $B_n$, $n = 1, 2, \ldots$, of compact open subgroups of $V \cap L$ with $B_1 = V \cap L$ and such that $B_n > B_{n+1}$ for each $n$ and $\bigcap_{n=1}^{\infty} B_n = \{e\}$. By [9, Lemma 1], it may further be supposed that each $B_n$ satisfies $\mathbf{T1}(\alpha|_L)$. Note that the restriction of $\alpha$ to $L$ is an automorphism.

To find an element of $\mathcal{K}$ in the coset $(V \cap L)x$ we construct an increasing sequence of integers $N_k$, $k = 1, 2, 3, \ldots$ such that the compact sets

(9)
$$\mathcal{C}_k := \{y \in (V \cap L)x \mid \text{for every } j \in \{1, \ldots, k\}, \, \alpha^n(y) \in B_j \text{ whenever } n \geq N_j\}$$

are not empty. The sequence is constructed inductively.

To begin, choose the coset representative $x$ to be in $\mathcal{L}$. That is possible because $L$ is the closure of $\mathcal{L}$. Then there is an integer $N_1$ such that $\alpha^n(x)$ is in $V \cap L$ whenever $n \geq N_1$. Thus, with this choice of $N_1$, $\mathcal{C}_1 \neq \emptyset$.

Next, suppose that $N_1, \ldots, N_k$ have been found so that $\mathcal{C}_k$ is not empty and choose $y$ in $\mathcal{C}_k$. Since $\alpha^n(y)$ is in $B_k$ for all $n \geq N_k$ the set, $\mathcal{A}$, of accumulation points of $\{\alpha^n(y)\}_{n \geq 0}$ is a non-empty closed $\alpha$-invariant subset of $B_k$. Choose an element, $a$, of $\mathcal{A}$. Then there is an $N_{k+1} > N_k$ such that $\alpha^{N_{k+1}}(y)$ belongs to $aB_{k+1}$. Put $y' = \alpha^{-N_{k+1}}(a^{-1})y$. Since $\mathcal{A}$ is $\alpha$-invariant and contained in $B_k$, it follows that $y'$ is in $(V \cap L)x$ and that, for $1 \leq j \leq k$, $\alpha^n(y')$ is in $B_j$ whenever $n \geq N_j$.

By construction, $\alpha^{N_{k+1}}(y')$ belongs to $B_{k+1}$. Since $B_{k+1}$ satisfies $\mathbf{T1}(\alpha|_L)$, $y'$ can be factored as

(10) $$y' = vy'',$$

where $\alpha^n(v)$ is in $B_{k+1}$ for every $n \leq N_{k+1}$ and $\alpha^n(y'')$ is in $B_{k+1}$ for every $n \geq N_{k+1}$. Then we have that $y''$ is in $(V \cap L)x$ because $y'' = v^{-1}y'$ with $v^{-1}$ in $B_{k+1}$ and $y'$ in $(V \cap L)x$. Furthermore, for each $j$ between 1 and $k+1$ we have that $\alpha^n(y'')$ is in $B_j$ whenever $n \geq N_j$ because: in case $N_j \leq n < N_{k+1}$, $\alpha^n(y'') = \alpha^n(v^{-1})\alpha^n(y')$ with $\alpha^n(v^{-1})$ in $B_{k+1}$ and $\alpha^n(y')$ in $B_j$; and, in case $n \geq N_{k+1}$, $\alpha^n(y'')$ is in $B_{k+1}$. Hence, with this choice of $N_{k+1}$, $y''$ belongs to $\mathcal{C}_{k+1}$ and this set is not empty. Therefore the inductive construction continues.

Since $\{\mathcal{C}_k\}_{k=1}^{\infty}$ is a decreasing sequence of non-empty compact sets, $\bigcap_{k=1}^{\infty} \mathcal{C}_k$ is not empty. Each element $y$ in this intersection belongs to $\mathcal{K}$ because membership of the intersection implies that $\lim_{j \to \infty} \alpha^j(y) = e$ and because $\{\alpha^{-j}(y)\}_{j \geq 0}$ is bounded for every $y$ in the compact and $\alpha$-invariant subgroup $L$. Therefore $L = (V \cap L)\mathcal{K}$ in the case when $L$ is separable.

When $L$ is not separable, choose $x$ in $\mathcal{L}$ as at the beginning of the recursive construction. Define $L'$ to be the closed subgroup generated by $\{\alpha^n(x)\}_{n \in \mathbf{Z}}$. Then $L'$ is separable and $\alpha$-invariant. Choose a sequence $B_1, B_2, \ldots$ of compact open subgroups of $L'$ satisfying $\mathbf{T1}(\alpha|_{L'})$ and with $B_1 = V \cap L'$. Then there is an $N_1$ such that $\alpha^n(x)$ is in $B_1$ for all $n \geq N_1$ and that is all that is needed to complete the recursive construction. □

**Lemma 2.4.** $W = V'L$ *is a maximal group contained in the set* $VL$.

**Proof.** Let $T$ be a group with $VL \supseteq T \supseteq V'L$. Then for each $t$ in $T$ we have that $t = vl$ for some $v$ in $V$ and $l$ in $L$. Since $L \leq T$ and $T$ is a group, $mvm^{-1}$ belongs to $T$ for each $m$ in $L$. Since $T \leq VL$, it follows that $v$ is in fact $V'$ and so $t$ belongs to $V'L$ as required. □

**Proof of Proposition 2.1.** Since $K$ is compact, [8, Lemma 3.3] shows that $V''K$ is a compact open subgroup of $G$. It must be shown that $V''K = V'L$.

Lemma 2.3 implies that $V'L = V'(V \cap L)K$, and this set equals $V'K$ because $V \cap L$ is contained in $V'$. The same reasoning shows that $VL = VK$ and it follows, by (3), (5) and Lemma 2.2, that $V' \leq V''$. Hence

$$V'L = V'K \leq V''K.$$

Now $V'' \leq V$ and, by Lemma 2.2, $K \leq L$. Hence $V''K$ is a group between $V'L$ and $VL$ and Lemma 2.4 shows that $V''K = W$. □

**Notation 2.5.** The present paper is concerned with compact open subgroups which are tidy for several automorphisms at once. Steps 1, 2a and 3a for producing a subgroup $W$ which is tidy for $\alpha$ will be called the *$\alpha$-tidying procedure*.

Should $U$ be tidy for automorphisms $\alpha$ and $\beta$ the factorings of $U$ given in **T1**$(\alpha)$ and **T1**$(\beta)$ may differ. To specify that the factoring given by $\alpha$, say, is being considered the notation

$$U_{\alpha+} = \bigcap_{k \geq 0} \alpha^k(U)$$

will be used, and similarly for $U_{\alpha-}$, $U_{\alpha 0} := U_{\alpha+} \cap U_{\alpha-}$, $U_{\alpha++}$ and $U_{\alpha--}$.

The following criteria for an element to belong to a tidy subgroup or to the contractive part of a tidy subgroup will be used repeatedly. They are implicit in [9] and [8] but are stated and proved explicitly here for ease of reference.

**Lemma 2.6.** *Let $W$ be a compact open subgroup of $G$ which is tidy for $\alpha$.*
  (i). *Let $w$ be in $W$ and suppose that $\{\alpha^n(w)\}_{n \geq 0}$ is bounded. Then $w$ belongs to $W_-$.*
  (ii). *Suppose that $\alpha^m(w)$ and $\alpha^n(w)$ belong to $W$ where $m < n$. Then $\alpha^k(w)$ is in $W$ for every $k$ between $m$ and $n$.*

**Proof.** (i) Since $W$ satisfies **T1**$(\alpha)$, $w = w_+ w_-$ where $w_\pm \in W_\pm$. Now $\{\alpha^n(w_-)\}_{n \geq 0}$ is contained in $W_-$, which is compact, and $w$ belongs to $W_-$ if and only if $w_+$ does, and so it may be supposed that $w$ belongs to $W_+$.

Let $c$ be an accumulation point of $\{\alpha^n(w)\}_{n \geq 0}$. Then $c$ belongs to the closure of $W_{++}$ and so, since $W$ satisfies **T2**$(\alpha)$, $c$ belongs to $W_{++}$. Hence there is a $k \geq 0$ such that $c$ is in $\alpha^k(W_+)$.

In the subspace topology on $W_{++}$, $\alpha^k(W_+)$ is an open neighbourhood of $c$. Hence there are arbitrarily large values of $n$ such that $\alpha^n(w)$ belongs to $\alpha^k(W_+)$. Applying $\alpha^{-k}$ yields that there are arbitrarily large values of $n$ such that $\alpha^{n-k}(w)$ belongs to $W_+$. However, if $\alpha^{n-k}(w)$ is in $W_+$, then $\alpha^j(w)$ belongs to it for every $j \leq n-k$. It follows that $\alpha^j(w)$ belongs to $W_+$ for every $j$ and hence that $w$ belongs to $W_-$.

(ii) By [9, Lemma 2] there are $u$ and $v$ in $G$ such that $w = uv$, $\alpha^m(u)$ belongs to $W_+$, $\alpha^n(u)$ belongs to $W_-$ and $\alpha^k(v)$ belongs to $W$ for $m \leq k \leq n$. It suffices then to show that $\alpha^k(u)$ is in $W$ for $m < k < n$.

Since $\alpha^n(u)$ is in $W_-$, $\{\alpha^k(u)\}_{k \geq n}$ is contained in $W_-$. Since this set is compact, $\{\alpha^k(u)\}_{k \geq m}$ is bounded and so, by part (i), $\alpha^m(u)$ belongs to $W_-$. Therefore $\alpha^k(u)$ belongs to $W$ for every $k \geq m$ and in particular for all $k$ between $m$ and $n$. □

## 3. Common Tidy Subgroups for Commuting Automorphisms

It is shown in this section that each finite set of mutually commuting elements of $\mathrm{Aut}(G)$ has a common tidy subgroup. The proof is by induction on the number of automorphisms. The inductive step is to show that, if $\alpha$ and $\beta$ are commuting elements of $\mathrm{Aut}(G)$ and $U$ is tidy for $\alpha$, then $\beta$-tidying $U$ produces a compact open subgroup $W$ which is tidy for both $\alpha$ and $\beta$.

The next result shows that Step 1 of $\beta$-tidying preserves $\alpha$-tidiness.

**Lemma 3.1.** *Let $\alpha$ and $\beta$ be commuting elements of $\mathrm{Aut}(G)$ and suppose that the compact open subgroup $U$ is tidy for $\alpha$. Then:*

(i). *$\beta^j(U)$ is tidy for $\alpha$ for each integer $j$ and*
(ii). *$\bigcap_{j=0}^n \beta^j(U)$ is tidy for $\alpha$ for each positive $n$.*

**Proof.** (i) Since $\alpha$ and $\beta$ commute, for each $j$ we have that

$$\begin{aligned}[\alpha(\beta^j(U)) : \alpha(\beta^j(U)) \cap \beta^j(U)] &= [\beta^j(\alpha(U)) : \beta^j(\alpha(U) \cap U)] \\ &= [\alpha(U) : \alpha(U) \cap U] \\ &= s(\alpha).\end{aligned}$$

Hence $\beta^j(U)$ is tidy for $\alpha$ by [8, Theorem 3.1]. Alternatively, it may be verified directly that $\beta^j(U)$ satisfies **T1**$(\alpha)$ and **T2**$(\alpha)$.

(ii) [9, Lemma 10] shows that the intersection of tidy subgroups is tidy. □

Next, we observe that the group $K$ introduced in Step 2a of $\beta$-tidying is $\alpha$-invariant. Recall also that Lemma 2.2 shows that this group is compact.

**Lemma 3.2.** *Let $\alpha$ and $\beta$ be commuting elements of $\mathrm{Aut}(G)$ and put*

$$K = \{g \in G : \lim_{j \to \infty} \beta^j(g) = e \text{ and } \{\beta^{-j}(g)\}_{j \geq 0} \text{ is bounded}\}^-.$$

*Then $K$ is a compact group and $\alpha(K) = K$.*

The $\alpha$-invariance of $K$ is needed in the proof that Step 3a preserves $\alpha$-tidiness. The group $L$ introduced in Step 2 of the earlier procedures is not easily seen to be $\alpha$-invariant and this is the reason for the new version. Several of the arguments in the following proof are based on those in [8]. Where possible results from [8] are simply cited.

**Theorem 3.3.** *Let $\alpha$ be in $\mathrm{Aut}(G)$ and suppose that the compact open subgroup $V$ is tidy for $\alpha$. Let $K$ be a compact subgroup of $G$ such that $\alpha(K) = K$ and put*

$$V'' = \{v \in V : lvl^{-1} \in VK \text{ for all } l \in K\}.$$

*Then $V''$ is an open subgroup of $G$.*

(i). *Let $v$ be in $V'' \cap V_+$ and $l$ be in $K$. Then $lvl^{-1}$ belongs to $V_+K$.*
(ii). *$V''_+ = V'' \cap V_+$*
(iii). *$V'' = V''_+ V''_-$.*
(iv). *$V''_{++} = \bigcup_{n \geq 0} \alpha^n(V''_+)$ is closed.*
(v). *$V''K$ is a compact open subgroup of $G$ which is tidy for $\alpha$.*

**Proof.** That $V''$ is an open subgroup of $G$ is shown in [8, Lemma 3.3].

(i) Since $v$ is in $V''$, $lvl^{-1} = uk$ for some $u$ in $V$ and $k$ in $K$. Then

$$\alpha^n(u) = \alpha^n(l)\alpha^n(v)\alpha^n(l^{-1}k^{-1}) \text{ for all } n$$

and so, since $K$ is $\alpha$-invariant and $\alpha^{-1}(V_+) \leq V_+$,

$$\alpha^n(u) \in KV_+K \text{ for all } n \leq 0.$$

Since $KV_+K$ is compact and $V$ is tidy for $\alpha$, it follows by Lemma 2.6(i) that $u$ belongs to $V_+$ as required.

(ii) The element $v$ belongs to $V''_+$ if and only if $\alpha^{-n}(v)$ is in $V''$ for every $n \geq 0$. Since $V'' \leq V$, this implies that $v$ belongs to $V_+$ and so

$$V''_+ \leq V'' \cap V_+.$$

Now let $v$ be in $V'' \cap V_+$. Let $l$ be in $K$ and $n \geq 0$. Then

$$l\alpha^{-n}(v)l^{-1} = \alpha^{-n}\left(\alpha^n(l)v\alpha^n(l)^{-1}\right)$$

where $\alpha^n(l)v\alpha^n(l)^{-1}$ belongs to $V_+K$ by part (i) and because $\alpha(K) = K$. Hence, using the $\alpha$-invariance of $K$ again,

$$l\alpha^{-n}(v)l^{-1} \in \alpha^{-n}(V_+)K \subset VK.$$

Therefore $\alpha^{-n}(v)$ belongs to $V''$ for every $n \geq 0$ and so $v$ belongs to $V''_+$.

(iii) Let $v$ be in $V''$. Then, since $V$ satisfies **T1**$(\alpha)$, there are $v_\pm \in V_\pm$ such that $v = v_+v_-$.

It may be supposed that there is $p \geq 0$ such that $\alpha^{-p}(v_+)$ belongs to $V''$. To see this, consider $\{\alpha^{-p}(v_+)\}_{p \geq 0}$. This set is contained in $V_+$ and so has an accumulation point, $x$ say. Then $\alpha^n(x)$ belongs to $V$ for every $n$ because it is an accumulation point of $\{\alpha^{n-p}(v_+)\}_{p \geq 0}$ and $\alpha^{n-p}(v_+)$ belongs to $V$ for all $p$ greater than $n$. Choose $p$ such that $\alpha^{-p}(v_+)$ is in $V''x$. Replace $v_+$ by $v_+\alpha^p(x)^{-1}$ and $v_-$ by $\alpha^p(x)v_-$. Then we still have that $v_\pm$ is in $V_\pm$ and $v = v_+v_-$ but now also have $\alpha^{-p}(v_+)$ in $V''$ as claimed.

It may now be shown that, under this supposition, $v_+$ belongs to $V''$. Let $l$ be in $K$. Then, since $\alpha^{-p}(v_+)$ is in $V''$ and $K$ is $\alpha$-invariant, part (i) shows that

(11) $$lv_+l^{-1} = um_1, \text{ where } u \in \alpha^p(V_+) \text{ and } m_1 \in K.$$

Hence

(12) $$lvl^{-1} = um_1(lv_-l^{-1}).$$

On the other hand, since $v$ is in $V''$, $lvl^{-1}$ belongs to $VK$ and so

(13) $$lvl^{-1} = w_+w_-m_2, \text{ for some } w_\pm \in V_\pm \text{ and } m_2 \in K.$$

Comparing (12) and (13) yields

$$w_+^{-1}u = w_-m_2(lv_-^{-1}l^{-1})m_1^{-1} \in U_-KU_-K.$$

It follows that $\{\alpha^n(w_+^{-1}u)\}_{n \geq 0}$ has an accumulation point but then, since $w_+^{-1}u$ belongs to $\alpha^p(V_+)$, Lemma 2.6(i) shows that $w_+^{-1}u$ is in $V_+ \cap V_-$. Since $w_+$ is in $V_+$, it follows that $u$ is in $V_+$. Hence, by (11),

$$lv_+l^{-1} \text{ belongs to } V_+K.$$

Since this holds for arbitrary $l$ in $K$, $v_+$ belongs to $V''$. It follows that $v_- = v_+^{-1}v$ is in $V''$ as well. Finally, (ii) implies that $v_\pm$ belongs to $V''_\pm$.

(iv) By [4, (5.37)] or [1, Prop. 2.4, Chapter III], $V''_{++}$ is closed if $V''_{++} \cap V''$ is closed. Hence it suffices to show that $V''_{++} \cap V'' = V''_+$. To this end, let $v$ be in $V''_+$ and suppose that $\alpha^n(v)$ belongs to $V''$ for some $n \geq 0$. Now $\alpha^n(v)$ being in $V''$, it may be factored as in part (iii),
$$\alpha^n(v) = w_+ w_- \text{ for some } w_\pm \in V''_\pm,$$
and it must be shown that $w_-$ belongs to $V''_+$.

Since $w_-$ belongs to $V''_-$, part (i) implies that for each $l$ in $K$

(14)  $$lw_- l^{-1} = uk \text{ for some } u \in V_- \text{ and } k \in K.$$

Now $\{\alpha^{-j}(u)\}_{j \geq 0}$ is bounded because
$$u = lw_- l^{-1} k^{-1} = lw_+^{-1} \alpha^n(v) l^{-1} k^{-1},$$
so that $\alpha^{-j}(u)$ belongs to the compact set $K\alpha^n(V''_+)K$ for each $j \geq 0$. Since $V$ is tidy for $\alpha$, it follows by Lemma 2.6(i) that $u$ belongs to $V_+ \cap V_-$. Therefore (14) shows that
$$lw_- l^{-1} \in (V_+ \cap V_-)K \text{ for every } l \in K.$$
Arguing as in (ii) we now see that $w_-$ belongs to $V''_+$ as required.

(v) That $V''K$ is a compact open subgroup is shown in [8, Lemma 3.3]. It is clear that $V''_+ \leq V''_+ K \leq (V''K)_+$ and similarly for $V''_-$. Since $V''$ satisfies **T1**$(\alpha)$ it follows that
$$V''K = V''_+ (V''_- K) = (V''K)_+ (V''K)_-$$
and so $V''K$ satisfies **T1**$(\alpha)$. Also, $\alpha^n(V''_+ K) = \alpha^n(V''_+)K$ for each $n$ because $K$ is $\alpha$-invariant and so
$$(V''K)_{++} = \bigcup_{n \geq 0} \alpha^n(V''_+ K) = \left( \bigcup_{n \geq 0} \alpha^n(V''_+) \right) K,$$
which is closed because $\bigcup_{n \geq 0} \alpha^n(V''_+)$ is closed and $K$ is compact. Hence $V''K$ satisfies **T2**$(\alpha)$. □

Now suppose that $\mathfrak{A}$ is a subset of $\mathrm{Aut}(G)$ such that there is a compact open subgroup $U$ which is tidy for each $\alpha$ in $\mathfrak{A}$. Let $\beta$ be an automorphism which commutes with each $\alpha$ in $\mathfrak{A}$. Then applying the $\beta$-tidying procedure to $U$ produces a subgroup $W$ which is tidy for $\beta$. Lemma 3.1 shows that the subgroup $V$ obtained in the first step of the procedure is tidy for each $\alpha$ in $\mathfrak{A}$ and Lemma 3.2 that the subgroup $K$ found in the second step is $\alpha$-invariant. Hence we are in a position to apply Theorem 3.3 in the third step and conclude that $W$ is tidy for each $\alpha$ in $\mathfrak{A}$ as well. This proves the following result.

**Theorem 3.4.** *Let $\mathfrak{A}$ be a set of automorphisms of the locally compact group $G$ and let $\beta$ be an automorphism which commutes with each $\alpha$ in $\mathfrak{A}$. Suppose that there is a compact open subgroup $U$ which is tidy for each $\alpha$ in $\mathfrak{A}$. Then there is a compact open subgroup $W$ tidy for each automorphism in $\mathfrak{A} \cup \{\beta\}$.*

An inductive argument now shows that each finite set of commuting automorphisms has a common tidy subgroup $U$. It does not follow however that this subgroup is tidy for every automorphism in the group of automorphisms generated by the finite set, as the following example shows.

**Example 3.5.** Let $F$ be a finite group with identity $e_F$ and let $\mathcal{I}$ be the indexing set $\mathbf{Z} \times (\mathbf{Z}/2\mathbf{Z})$, where $\mathbf{Z}/2\mathbf{Z} = \{\bar{0}, \bar{1}\}$. Define

$$G = \{f \in F^{\mathcal{I}} : \exists N \in \mathbf{Z} \text{ such that } f(n, \bar{a}) = e_F \text{ whenever } n < N\}.$$

Then $G$ is a group under coordinatewise multiplication. For each $n$ define the subgroup $G_n = \{f \in G : f(m, \bar{a}) = e_F \text{ for all } m < n\}$. Define a topology on $G$ by letting $\{G_n : n \in \mathbf{Z}\}$ be a base of neigbourhoods for $e_G$. Then $G$ is a totally disconnected locally compact group.

Let $\alpha_1$ and $\alpha_2$ be the automorphisms

$$\alpha_1(f)(n, \bar{a}) = f(n+1, \bar{a}), \qquad \alpha_2(f)(n, \bar{a}) = f(n+1, \bar{a} + \bar{1}).$$

Then $\alpha_1$ and $\alpha_2$ commute. Let $U$ be the compact open subgroup

$$U = \{f \in G : f(n, \bar{a}) = e_F \text{ if } n < 0 \text{ or if } n = 0 \text{ and } \bar{a} = \bar{0}\}.$$

Then

$$\alpha_1(U) = \{f \in G : f(n, \bar{a}) = e_F \text{ if } n < -1 \text{ or if } n = -1 \text{ and } \bar{a} = \bar{0}\}$$

and

$$\alpha_2(U) = \{f \in G : f(n, \bar{a}) = e_F \text{ if } n < -1 \text{ or if } n = -1 \text{ and } \bar{a} = \bar{1}\}.$$

Hence $\alpha_1(U) > U$ and $\alpha_2(U) > U$ so that $U$ is tidy for both $\alpha_1$ and $\alpha_2$.

On the other hand,

$$\alpha_1^{-1}\alpha_2(f)(n, \bar{a}) = f(n, \bar{a} + \bar{1}).$$

It follows that

$$U \cap \alpha_1^{-1}\alpha_2(U) = \{f \in G : f(n, \bar{a}) = e_F \text{ if } n \leq 0\}$$

and that this group is invariant under $\alpha_1^{-1}\alpha_2$. Therefore $U$ is not tidy for $\alpha_1^{-1}\alpha_2$.

The above example notwithstanding, it will be shown in §5 that each finitely generated abelian group of automorphisms has a common tidy subgroup. The proof relies on a refinement of property **T1** which will now be described.

## 4. Factoring Subgroups Tidy for a Finite Set of Automorphisms

Property **T1**$(\alpha)$ expresses a compact open subgroup $U$ which is tidy for $\alpha$ as the product of a subgroup on which $\alpha$ is expanding and a subgroup on which $\alpha$ is contracting and, if $U$ is tidy for several automorphisms, there are several such factorings of $U$. In this section it is shown that, under certain conditions on the automorphisms given in Definition 4.2, $U$ may be written as the product of subgroups such that each automorphism is either expanding or contracting on each subgroup.

The conditions in 4.2 are satisfied if the automorphisms commute. In the converse direction, it is seen in this section, if $\mathfrak{H}$ is a group of automorphisms satisfying the conditions in 4.2, then $\mathfrak{H}$ is abelian modulo those automorphisms which leave tidy subgroups invariant.

As is the case for a single automorphism, the factors of the tidy subgroup $U$ correspnding to several automorphisms are intersections of images of $U$ under powers of the automorphisms.

**Definition 4.1.** For each finite set $\mathfrak{a}$ of automorphisms of $G$ and each compact open subgroup $U$ of $G$ define
$$U_{\mathfrak{a}} := \bigcap_{\alpha \in \mathfrak{a}} U_{\alpha+}.$$
If $\beta$ is another automorphism, put
$$U_{\mathfrak{a},\beta+} := (U_{\mathfrak{a}})_{\beta+} = \bigcap_{k \geq 0} \beta^k(U_{\mathfrak{a}}),$$
and define $U_{\mathfrak{a},\beta-}$ and $U_{\mathfrak{a},\beta 0}$ similarly.

In the following, we shall need to know that $\alpha(U_{\mathfrak{a}}) \geq U_{\mathfrak{a}}$ for each $\alpha$ in $\mathfrak{a}$ but that is not true of general finite sets of automorphisms. It will be shown in Theorem 4.6 however that it does hold for sets of automorphisms satisfying the condition given next. The condition extends to finite sets of automorphisms a fact which is easily proved for a single automorphism $\alpha$, namely, that if $U$ is tidy for $\alpha$, then $\alpha^n(U)$ is tidy for $\alpha$ for each $n$.

**Definition 4.2.** Let $U$ be a compact open subgroup of $G$ and $\mathfrak{H}$ be a group of automorphisms of $G$.
  (i). $U$ is *invariably tidy* for the finite set, $\mathfrak{a}$, of automorphisms means that for every $\gamma$ in $\langle \mathfrak{a} \rangle$, $\gamma(U)$ is tidy for each $\alpha$ in $\mathfrak{a}$.
  (ii). $\mathfrak{H}$ is said to have *local tidy subgroups* if for every finite subset, $\mathfrak{a}$, of $\mathfrak{H}$ there is a compact open subgroup $U$ which is invariably tidy for $\mathfrak{a}$.
  (iii). $\mathfrak{H}$ is said to have *tidy subgroups* if there is a compact open subgroup $U$ which is tidy for every $\alpha$ in $\mathfrak{H}$. In this case $U$ is said to be *tidy for $\mathfrak{H}$*.

Theorem 3.4 and Lemma 3.1 imply that there is, for each finite set $\mathfrak{a}$ of commuting automorphisms, a compact open subgroup which is invariably tidy for $\mathfrak{a}$. Hence abelian groups of automorphisms have local tidy subgroups. As the terminology suggests, the property of having tidy subgroups is stronger than having local tidy subgroups.

**Lemma 4.3.** *Let $\alpha$ be an element of $\mathfrak{H}$ and suppose that $U$ is tidy for $\mathfrak{H}$. Then $\alpha(U)$ is tidy for $\mathfrak{H}$.*

**Proof.** Let $\beta$ be any automorphism in $\mathrm{Aut}(G)$. Then $\alpha(U)$ is tidy for $\alpha\beta\alpha^{-1}$ as may be seen by verifying properties $\mathbf{T1}(\alpha\beta\alpha^{-1})$ and $\mathbf{T2}(\alpha\beta\alpha^{-1})$ directly or by observing that $\alpha(U)$ minimises the index $[\alpha\beta\alpha^{-1}(V) : V \cap \alpha\beta\alpha^{-1}(V)]^1$. Since $U$ is tidy for $\mathfrak{H}$, it follows that $\alpha(U)$ is tidy for $\alpha\mathfrak{H}\alpha^{-1} = \mathfrak{H}$. □

Let $G$, $\alpha_1$, $\alpha_2$ and $U$ be as in Example 3.5. Then Lemma 3.1 shows that $U$ is invariably tidy for $\mathfrak{a} = \{\alpha_1, \alpha_2\}$. However $U$ is not tidy for $\langle \alpha_1, \alpha_2 \rangle$ and so it is not immediately obvious that having local tidy subgroups should imply that $\mathfrak{H}$ has tidy subgroups. It will be shown in §5 though that, nevertheless, this is the case if $\mathfrak{H}$ is finitely generated.

**Lemma 4.4.** *Let $\mathfrak{a}$ be a finite set of automorphisms of $G$. Let $\beta$ be another automorphism and $U$ be a compact open subgroup of $G$ such that $\alpha^k(U)$ is tidy for $\beta$ for each $\alpha \in \mathfrak{a}$ and each $k \geq 0$. Then:*
  (i). $U_{\mathfrak{a},\beta+} = U_{\mathfrak{a} \cup \{\beta\}}$; *and*

---
[1] I am grateful to Harald Biller for this simplified argument.

(ii). $U_{\mathfrak{a}} = U_{\mathfrak{a},\beta-}U_{\mathfrak{a},\beta+}$.

**Proof.** For each $m \geq 0$, put $U_{\mathfrak{a}}^{(m)} = \bigcap_{\alpha \in \mathfrak{a}} \bigcap_{k=0}^{m} \alpha^k(U)$. Then $U_{\mathfrak{a}}^{(m)}$ is a compact open subgroup of $G$ and, since $\alpha^k(U)$ is tidy for $\beta$ for each $\alpha$ and $k$, [9, Lemma 10] shows that $U_{\mathfrak{a}}^{(m)}$ is tidy for $\beta$. We have that $U_{\mathfrak{a}} = \bigcap_{m \geq 0} U_{\mathfrak{a}}^{(m)}$.

(i) Since $U_{\mathfrak{a}} \leq U$, it is clear that $U_{\mathfrak{a},\beta+} \leq U_{\beta+}$ and that $U_{\mathfrak{a},\beta+} \leq U_{\mathfrak{a}}$ follows immediately from the definition. Hence $U_{\mathfrak{a},\beta+} \leq U_{\mathfrak{a}\cup\{\beta\}}$.

Now let $u$ be in $U_{\mathfrak{a}\cup\{\beta\}}$. Then $\{\beta^{-k}(u)\}_{k\geq 0}$ is bounded because $u$ belongs to $U_{\beta+}$. Since $U_{\mathfrak{a}}^{(m)}$ is tidy for $\beta$, Lemma 2.6(i) shows that $u$ belongs to $\left(U_{\mathfrak{a}}^{(m)}\right)_{\beta+}$ for each $m$. Therefore $u$ is in $U_{\mathfrak{a},\beta+}$.

(ii) Since $U_{\mathfrak{a}}^{(m)}$ is tidy for $\beta$, we have

$$
(15) \qquad U_{\mathfrak{a}}^{(m)} = \left(U_{\mathfrak{a}}^{(m)}\right)_{\beta-}\left(U_{\mathfrak{a}}^{(m)}\right)_{\beta+}.
$$

Let $u$ be in $U_{\mathfrak{a}}$ and for each $m \geq 0$ put

$$
C_m = \left\{(u_-, u_+) \in U_{\mathfrak{a}}^{(m)} \times U_{\mathfrak{a}}^{(m)} : u_\pm \in \left(U_{\mathfrak{a}}^{(m)}\right)_{\beta\pm} \text{ and } u = u_-u_+\right\}.
$$

Then, since $U_{\mathfrak{a}} \leq U_{\mathfrak{a}}^{(m)}$, (15) shows that $C_m$ is non-empty for each $m$. Now $C_m$ is compact and $C_{m+1} \subseteq C_m$ for each $m$ and so $\bigcap_{m \geq 0} C_m \neq \emptyset$. Choose $(u_-, u_+)$ in this intersection. Then $u_+$ and $u_-$ belong to $\bigcap_{m \geq 0} U_{\mathfrak{a},\beta\pm}^{(m)} = U_{\mathfrak{a},\beta\pm}$ and $u = u_-u_+$. $\square$

Note that the Lemma does not suppose that $U$ is tidy for the automorphisms in $\mathfrak{a}$. This observation will be needed later on.

Some of the arguments to follow will consider a sequence of automorphisms having an invariably tidy subgroup $U$ and part (ii) of the lemma will be used repeatedly to factor $U$. Here is notation which will be used in those arguments.

**Definition 4.5.** Given the sequence $\alpha_1, \ldots, \alpha_n$ of automorphisms of $G$ define, for each function $\epsilon : \{1, \ldots, n\} \to \{-1, 1\}$, the set of automorphisms

$$
\mathfrak{a}_\epsilon := \left\{\alpha_j^{\epsilon(j)} : j \in \{1, \ldots, n\}\right\}.
$$

The collection of these sets of automorphisms is denoted

$$
\Phi := \left\{\mathfrak{a}_\epsilon : \epsilon \in \{-1, 1\}^{\{1, \ldots, n\}}\right\}.
$$

The lexicographic order on the indexing set $\{-1, 1\}^{\{1, \ldots, n\}}$ (where $-1 < 1$ as usual) determines a total order on $\Phi$. This will be called the *halving order* on $\Phi$.

Now suppose that $U$ is invariably tidy for the set of automorphisms $\{\alpha_1, \ldots, \alpha_n\}$. Then $U = U_{\alpha_1-}U_{\alpha_1+}$ because $U$ is tidy for $\alpha_1$. Since $U$ is invariably tidy for $\{\alpha_1, \ldots, \alpha_n\}$, Lemma 4.4(ii) may be applied repeatedly with $\beta$ equal in turn to $\alpha_2$, $\ldots$, $\alpha_n$ until $U$ is written as a product of $2^n$ factors. Lemma 4.4(i) shows that each of the factors has the form $U_{\mathfrak{a}}$ for some $\mathfrak{a}$ in $\Phi$. Moreover, the order of these $2^n$ factors is the halving order on $\Phi$. Since $\beta$ is expanding on $U_\beta$, we have that each $\alpha_j$ is expanding or contracting on $U_{\mathfrak{a}_\epsilon}$ depending on whether $\epsilon(j)$ equals 1 or $-1$. This procedure thus yields the following decomposition of invariably tidy subgroups.

**Theorem 4.6.** *Let $\alpha_1, \ldots, \alpha_n$ be a sequence of automorphisms of $G$ and suppose that $U$ is a compact open subgroup of $G$ which is invariably tidy for $\{\alpha_1, \ldots, \alpha_n\}$. Then*

(16) $$U = \prod_{\mathfrak{a} \in \Phi} U_\mathfrak{a},$$

*where the factors $U_\mathfrak{a}$ in the product are in the halving order. For each $j \in \{1, \ldots, n\}$ we have $\alpha_j(U_\mathfrak{a}) \geq U_\mathfrak{a}$ or $\alpha_j(U_\mathfrak{a}) \leq U_\mathfrak{a}$ depending on whether $\alpha_j \in \mathfrak{a}$ or $\alpha_j^{-1} \in \mathfrak{a}$. The equation* 16 *will be called the* factoring of $U$ corresponding to the sequence of automorphisms $\alpha_1, \ldots, \alpha_n$.

If $U$ is invariably tidy for $\{\alpha_1, \ldots \alpha_n\}$, then $\beta(U)$ is invariably tidy for $\{\alpha_1, \ldots \alpha_n\}$ for every $\beta$ in $\langle \alpha_1, \ldots \alpha_n \rangle$. Hence $\beta(U)$ also has a factoring, $\beta(U) = \prod_{\mathfrak{a} \in \Phi} \beta(U)_\mathfrak{a}$, corresponding to $\alpha_1, \ldots, \alpha_n$. The next few results, leading to Corollary 4.9, show that this factoring is obtained by applying $\beta$ to (16).

**Lemma 4.7.** *Let the compact open subgroup $U$ be invariably tidy for $\{\alpha, \beta\}$. Then*

(17) $$\beta\left(U_{\alpha+} \cap U_{\beta-}\right) \leq \beta(U)_{\alpha+} \cap \beta(U)_{\beta-}$$

**Proof.** Let $u$ be an element of $\beta\left(U_{\alpha+} \cap U_{\beta-}\right)$. Then $\beta(u)$ is in $\beta(U)_{\beta-}$ (in fact in $U_{\beta-}$) because $u$ is in $U_{\beta-}$. To see that $\beta(u)$ belongs to $\beta(U)_{\alpha+}$ as well, write $U_{\alpha+}$ as the intersection of subgroups

$$U_{\alpha+} = \bigcap_{m \geq 0} U_{\{\alpha\}}^{(m)}, \quad \text{where } U_{\{\alpha\}}^{(m)} = \bigcap_{k=0}^{m} \alpha^k(U).$$

Each $U_{\{\alpha\}}^{(m)}$ is tidy for $\beta$ by Lemma 10 in [9] and $u$ belongs to each subgroup. The sequence $\{\beta^n(u)\}_{n \geq 0}$ is relatively compact because $u$ is in $U_{\beta-}$ and so, by Lemma 2.6(i), $u$ belongs to $\left(U_{\{\alpha\}}^{(m)}\right)_{\beta-}$ for each $m$. Hence $\beta(u)$ belongs to $U_{\{\alpha\}}^{(m)}$ for each $m$ and it follows that $\beta(u)$ is in $U_{\alpha+}$. Therefore the sequence $\{\alpha^{-n}\beta(u)\}_{n \geq 0}$ is relatively compact and we conclude with the aid of Lemma 2.6(i) that $\beta(u)$ belongs to $\beta(U)_{\alpha+}$. □

**Theorem 4.8.** *Let the compact open subgroup $U$ be invariably tidy for $\{\alpha, \beta\}$. Then*

$$\beta(U_{\alpha+}) = \beta(U)_{\alpha+}.$$

**Proof.** First observe that it suffices to show that

(18) $$\beta(U_{\alpha+}) \leq \beta(U)_{\alpha+}.$$

For this, recall that $\gamma(U)$ is tidy for $\beta$ if and only if it is tidy for $\beta^{-1}$. Hence $U$ is in fact invariably tidy for $\{\alpha, \beta, \alpha^{-1}, \beta^{-1}\}$. Now taking $\beta(U)$ in place of $U$ and $\beta^{-1}$ in place of $\beta$ in (18) yields $\beta^{-1}(\beta(U)_{\alpha+}) \leq U_{\alpha+}$. Applying $\beta$ to each side shows the reverse inclusion to (18).

Since $U$ is invariably tidy for $\{\alpha, \beta\}$, Lemma 4.4 shows that

$$\beta(U_{\alpha+}) = \beta(U_{\alpha+} \cap U_{\beta-}) \beta(U_{\alpha+} \cap U_{\beta+}) \quad \text{and}$$

$$\beta(U)_{\alpha+} = \left(\beta(U)_{\alpha+} \cap \beta(U)_{\beta-}\right)\left(\beta(U)_{\alpha+} \cap \beta(U)_{\beta+}\right).$$

Hence, by Lemma 4.7, (18) will follow once it is shown that

(19) $$\beta\left(U_{\alpha+} \cap U_{\beta+}\right) \leq \beta(U)_{\alpha+} \cap \beta(U)_{\beta+}.$$

To show this, let $u$ be in $U_{\alpha+} \cap U_{\beta+}$. Then $\beta(u)$ belongs to $\beta(U)_{\beta+}$ and so, by Lemma 4.4,

(20) $$\beta(u) = v_+ v_-, \text{ where } v_\pm \in \beta(U)_{\alpha\pm} \cap \beta(U)_{\beta+}.$$

We claim that

(21) $$\alpha^{-l}(v_-) \text{ belongs to } U\beta(U) \text{ for every } l \geq 0.$$

Towards establishing the claim, let $l$ be a positive integer. Since $v_-$ is in $\beta(U)_{\alpha-}$, $\{\alpha^m(v_-)\}_{m\geq 0}$ has an accumulation point, $t$ say, and $t$ belongs to $\beta(U)_{\alpha+} \cap \beta(U)_{\alpha-}$. Now $\alpha^m(v_-)$ also belongs to $\alpha^m \beta(U)_{\beta+}$ for every $m \geq 0$, by Lemma 4.7 and is clearly in $\beta(U)$. Since $\beta(U)$ is tidy for $\beta$, it follows that $\alpha^m(v_-)$ is in $\beta(U)_{\beta+}$ for every $m \geq 0$. Hence

(22) $$t \text{ belongs to } \beta(U)_{\alpha+} \cap \beta(U)_{\beta+}.$$

Now choose an $m \geq 0$ such that

(23) $$\alpha^m(v_-)t^{-1} \text{ belongs to the open subgroup } \beta^{-1}\alpha^l(U).$$

Lemma 4.7 and (22) imply that $\alpha^{-m}(t)$ belongs to $\alpha^{-m}\beta(U)_{\alpha+} \cap \alpha^{-m}\beta(U)_{\beta+}$ and it is clear that it also belongs to $\beta(U)$. Since $\beta(U)$ is tidy for $\alpha$ and $\beta$, it follows that $\alpha^{-m}(t)$ belongs to $\beta(U)_{\alpha+} \cap \beta(U)_{\beta+}$. Hence, by Lemma 4.7 again, $\beta^{-1}\alpha^{-m}(t)$ belongs to $U_{\alpha+} \cap U_{\beta+}$. Now $\beta^{-1}(v_-)$ also belongs to $U_{\alpha+} \cap U_{\beta+}$ because $\beta^{-1}(v_-) = \beta^{-1}(v_+)u$, by (20), where $u$ was chosen from $U_{\alpha+} \cap U_{\beta+}$ and $\beta^{-1}(v_+)$ is in $U_{\alpha+} \cap U_{\beta+}$ by (20) and Lemma 4.7. Therefore

(24) $$\beta^{-1}\alpha^{-m}\left(\alpha^m(v_-)t^{-1}\right) \text{ belongs to } U_{\alpha+} \cap U_{\beta+}.$$

Put $w = \alpha^{-(l+m)}\left(\beta^{-1}\alpha^{-m}\left(\alpha^m(v_-)t^{-1}\right)\right)$. Then (24) implies that $w$ belongs to $U$. We also have that
$$\left(\alpha^{-l}\beta\alpha^{l+m}\right)^2(w) = \alpha^{-l}\beta\left(\alpha^m(v_-)t^{-1}\right),$$
which belongs to $U$ by (23). Since $U$ is tidy for $\alpha^{-l}\beta\alpha^{l+m}$, it follows by Lemma 2.6(ii) that $\alpha^{-l}\beta\alpha^{l+m}(w)$ belongs to $U$. However
$$\alpha^{-l}\beta\alpha^{l+m}(w) = \alpha^{-l}(v_-)\alpha^{-(l+m)}(t^{-1})$$
and $\alpha^{-(l+m)}(t^{-1})$ is in $\beta(U)$ by (22). The claim (21) is thus established.

It follows from (21) that $\{\alpha^{-l}(v_-)\}_{l\geq 0}$ is relatively compact and so, since $\beta(U)$ is tidy for $\alpha$, Lemma 2.6(i) shows that $v_-$ is in $\beta(U)_{\alpha+}$. By construction, $v_+$ is in $\beta(U)_{\alpha+}$ too and so $\beta(u)$ belongs to this subgroup. That $\beta(u)$ belongs to $\beta(U)_{\beta+}$ is clear. Therefore (19) holds and the proof of the Theorem is complete. □

By definition, $U_\mathfrak{a} = \bigcap_{\alpha \in \mathfrak{a}} U_{\alpha+}$ for any finite set, $\mathfrak{a}$, of automorphisms. An immediate consequence of the Theorem then is the corresponding statement for finite sets of automorphisms.

**Corollary 4.9.** *Let the compact open subgroup $U$ be invariably tidy for the finite set $\mathfrak{a} \cup \{\beta\}$. Then*
$$\beta(U_\mathfrak{a}) = \beta(U)_\mathfrak{a}.$$

If $\mathfrak{a} \cup \{\beta\}$ is a commuting set of automorphisms, then Lemma 3.1 shows that any $U$ which is tidy for $\mathfrak{a} \cup \{\beta\}$ is invariably tidy. The corollary thus applies to commuting sets of automorphisms. However, when the automorphisms commute it is not difficult to show directly that $\beta(U_\mathfrak{a}) = \beta(U)_\mathfrak{a}$.

If $U$ is invariant under each automorphism in $\mathfrak{a} \cup \{\beta\}$, that is, if $\alpha(U) = U$ for each $\alpha$ in $\mathfrak{a} \cup \{\beta\}$, then $U$ is invariably tidy for $\mathfrak{a} \cup \{\beta\}$. In this case the conclusion of the corollary is obvious too.

Theorem 4.14 below shows that sets of automorphisms having an invariably tidy subgroup are commuting modulo automorphisms which leave $U$ invariant. The proof uses Corollary 4.9 and the following observations about the subgroup $U_\mathfrak{a}$.

**Lemma 4.10.** *Let the compact open group $U$ be invariably tidy for the finite set of automorphisms $\mathfrak{a} \cup \{\beta\}$. Then:*

(i). *$\beta(U_\mathfrak{a}) \cap U_\mathfrak{a} = \beta(U_\mathfrak{a}) \cap U$ and is an open subgroup of $\beta(U_\mathfrak{a})$;*
(ii). *$\beta(U_\mathfrak{a}) \cap U_\mathfrak{a} = \beta(U_{\mathfrak{a},\beta-})U_{\mathfrak{a},\beta+}$;*
(iii). *$[\beta(U_\mathfrak{a}) : \beta(U_\mathfrak{a}) \cap U_\mathfrak{a}] = [\beta(U_{\mathfrak{a},\beta+}) : U_{\mathfrak{a},\beta+}]$; and*
(iv). *$U_{\mathfrak{a},\beta++}$ is a closed subgroup of $G$.*

**Proof.** (i) By Corollary 4.9, $\beta(U_\mathfrak{a}) \cap U = \beta(U)_\mathfrak{a} \cap U$ and $\beta(U)_\mathfrak{a} \cap U = \beta(U)_\mathfrak{a} \cap U_\mathfrak{a}$, by Lemma 2.6(i). Hence $\beta(U_\mathfrak{a}) \cap U \leq \beta(U_\mathfrak{a}) \cap U_\mathfrak{a}$. The reverse inclusion is obvious.

(ii) Lemma 4.4(ii) shows that $\beta(U_\mathfrak{a}) = \beta(U_{\mathfrak{a},\beta-})\beta(U_{\mathfrak{a},\beta+})$. Since $\beta$ is expanding on $U_{\mathfrak{a},\beta+}$, it follows that $\beta(U_\mathfrak{a}) \geq \beta(U_{\mathfrak{a},\beta-})U_{\mathfrak{a},\beta+}$. Since $\beta$ is contracting on $U_{\mathfrak{a},\beta-}$, $\beta(U_{\mathfrak{a},\beta-}) \leq U_{\mathfrak{a},\beta-}$ and so $U_\mathfrak{a} \geq \beta(U_{\mathfrak{a},\beta-})U_{\mathfrak{a},\beta+}$. Hence $\beta(U_\mathfrak{a}) \cap U_\mathfrak{a} \geq \beta(U_{\mathfrak{a},\beta-})U_{\mathfrak{a},\beta+}$.

For the reverse inclusion, let $u$ be in $\beta(U_\mathfrak{a}) \cap U_\mathfrak{a}$. Then, since $u$ is in $\beta(U_\mathfrak{a})$,

$$u = u_- u_+, \quad \text{where } u_- \in \beta(U_{\mathfrak{a},\beta-}) \leq U_{\mathfrak{a},\beta-} \text{ and } u_+ \in \beta(U_{\mathfrak{a},\beta+}).$$

Hence $u_+ = u_-^{-1} u$ belongs to $\beta(U_{\mathfrak{a},\beta+}) \cap U_\mathfrak{a}$. Appealing to Corollary 4.9 and Lemma 2.6(i) again, we see that $u_+$ belongs to $U_{\mathfrak{a},\beta+}$.

(iii) Define a map $\Theta : \beta(U_{\mathfrak{a},\beta+})/U_{\mathfrak{a},\beta+} \to \beta(U_\mathfrak{a})/(\beta(U_\mathfrak{a}) \cap U_\mathfrak{a})$ by

$$\Theta(\beta(u)U_{\mathfrak{a},\beta+}) = \beta(u)(\beta(U_\mathfrak{a}) \cap U_\mathfrak{a}).$$

Then $\Theta$ is well-defined because $U_{\mathfrak{a},\beta+} \leq \beta(U_\mathfrak{a}) \cap U_\mathfrak{a}$ and is surjective because $\beta(U_\mathfrak{a}) = \beta(U_{\mathfrak{a},\beta+})\beta(U_{\mathfrak{a},\beta-})$ and $\beta(U_{\mathfrak{a},\beta-}) \leq \beta(U_\mathfrak{a}) \cap U_\mathfrak{a}$. That $\Theta$ is injective follows because $\beta(U_{\mathfrak{a},\beta+}) \cap \beta(U_\mathfrak{a}) \cap U_\mathfrak{a} = U_{\mathfrak{a},\beta+}$, which is in turn a consequence of Corollary 4.9 and Lemma 2.6(i).

(iv) By definition $U_{\mathfrak{a},\beta++} = \bigcup_{l \geq 0} \beta^l(U_{\mathfrak{a},\beta+})$ and it follows, by Corollary 4.9 and Lemma 2.6(i), that

$$U_{\mathfrak{a},\beta++} = \bigcup_{l \geq 0} \left\{ u \in \beta^l(U_{\beta+}) : \{\alpha^{-k}(u)\}_{k \geq 0} \text{ is bounded for every } \alpha \in \mathfrak{a} \right\}$$

$$= \left\{ u \in U_{\beta++} : \{\alpha^{-k}(u)\}_{k \geq 0} \text{ is bounded for every } \alpha \in \mathfrak{a} \right\}.$$

Now for each $\alpha$, $\left\{ u \in G : \{\alpha^{-k}(u)\}_{k \geq 0} \text{ is bounded} \right\}$ is a closed subgroup of $G$, see [9, Proposition 3], and $U_{\beta++}$ is also closed. Hence $U_{\mathfrak{a},\beta++}$ is the intersection of closed subgroups and is closed. □

**Lemma 4.11.** *Let $\mathfrak{a}$ be a finite set of automorphisms of $G$ and let $U$ and $V$ be compact open subgroups tidy for $\mathfrak{a}$. Suppose that $V \leq U$. Then $V_\mathfrak{a} = U_\mathfrak{a} \cap V$ and the index of $V_\mathfrak{a}$ in $U_\mathfrak{a}$ is finite.*

**Proof.** It is immediate from the definitions that $V_{\mathfrak{a}} \leq U_{\mathfrak{a}} \cap V$. For the reverse inclusion, let $u$ be in $U_{\mathfrak{a}} \cap V$. Then $\{\alpha^{-m}(u)\}_{m\geq 0}$ is bounded for every $\alpha$ in $\mathfrak{a}$. Since $u$ is also in $V$ and $V$ is tidy for each $\alpha$ in $\mathfrak{a}$, Lemma 2.6(i) shows that $u$ belongs to $V_{\mathfrak{a}}$.

It follows from the first part that $V_{\mathfrak{a}}$ is an open subgroup of the compact group $U_{\mathfrak{a}}$. Hence $[U_{\mathfrak{a}} : V_{\mathfrak{a}}]$ is finite. □

**Theorem 4.12.** *Let the compact open subgroups $U$ and $V$ be invariably tidy for the finite set $\mathfrak{a} \cup \{\beta\}$ of automorphisms. Then*

$$[\beta(V_{\mathfrak{a},\beta+}) : V_{\mathfrak{a},\beta+}] = [\beta(U_{\mathfrak{a},\beta+}) : U_{\mathfrak{a},\beta+}] < \infty.$$

**Proof.** Lemma 4.10(i) implies that $U_{\mathfrak{a},\beta+}$ is an open subgroup of $\beta(U_{\mathfrak{a},\beta+})$. Hence $[\beta(U_{\mathfrak{a},\beta+}) : U_{\mathfrak{a},\beta+}]$ is finite.

To show that $[\beta(U_{\mathfrak{a},\beta+}) : U_{\mathfrak{a},\beta+}]$ is equal to $[\beta(V_{\mathfrak{a},\beta+}) : V_{\mathfrak{a},\beta+}]$ it suffices to show that both equal $[\beta((U \cap V)_{\mathfrak{a},\beta+}) : (U \cap V)_{\mathfrak{a},\beta+}]$. Since the intersection of tidy subgroups is tidy, $U \cap V$ is invariably tidy for $\mathfrak{a} \cup \{\beta\}$. Hence it may be supposed, without loss of generality, that $V \leq U$.

Then $V_{\mathfrak{a},\beta+} \leq \beta(V_{\mathfrak{a},\beta+}) \leq \beta(U_{\mathfrak{a},\beta+})$ and it follows that

(25) $\qquad [\beta(U_{\mathfrak{a},\beta+}) : V_{\mathfrak{a},\beta+}] = [\beta(U_{\mathfrak{a},\beta+}) : \beta(V_{\mathfrak{a},\beta+})][\beta(V_{\mathfrak{a},\beta+}) : V_{\mathfrak{a},\beta+}].$

Similarly, $V_{\mathfrak{a},\beta+} \leq U_{\mathfrak{a},\beta+} \leq \beta(U_{\mathfrak{a},\beta+})$ and

(26) $\qquad [\beta(U_{\mathfrak{a},\beta+}) : V_{\mathfrak{a},\beta+}] = [\beta(U_{\mathfrak{a},\beta+}) : U_{\mathfrak{a},\beta+}][U_{\mathfrak{a},\beta+} : V_{\mathfrak{a},\beta+}].$

Now $[U_{\mathfrak{a},\beta+} : V_{\mathfrak{a},\beta+}]$ is finite by Lemma 4.11 and is equal to $[\beta(U_{\mathfrak{a},\beta+}) : \beta(V_{\mathfrak{a},\beta+})]$ because $\beta$ is an isomorphism. Therefore (25) and (26) imply that

$$[\beta(V_{\mathfrak{a},\beta+}) : V_{\mathfrak{a},\beta+}] = [\beta(U_{\mathfrak{a},\beta+}) : U_{\mathfrak{a},\beta+}].$$

□

**Definition 4.13.** The *scale of $\beta$ relative to $\mathfrak{a}$* is the positive integer

$$s_{\mathfrak{a}}(\beta) := [\beta(U_{\mathfrak{a},\beta+}) : U_{\mathfrak{a},\beta+}] = [\beta(U_{\mathfrak{a}}) : \beta(U_{\mathfrak{a}}) \cap U_{\mathfrak{a}}],$$

where $U$ is a compact open subgroup invariably tidy for $\mathfrak{a} \cup \{\beta\}$.

It may now be shown that automorphisms having an invariably tidy subgroup commute modulo the normaliser of that subgroup.

**Theorem 4.14.** *Let $\alpha$ and $\beta$ be automorphisms of $G$ and let the compact open subgroup $U$ be invariably tidy for $\{\alpha, \beta, \alpha\beta\alpha^{-1}\beta^{-1}\}$. Then $\alpha\beta\alpha^{-1}\beta^{-1}(U) = U$.*

**Proof.** Put $\alpha_1 = \alpha$, $\alpha_2 = \beta$ and $\alpha_3 = \alpha\beta\alpha^{-1}\beta^{-1}$ and let

(27) $$U = \prod_{\mathfrak{a} \in \Phi} U_{\mathfrak{a}},$$

be the corresponding factoring of $U$. Recall that $\Phi = \{\mathfrak{a}_\epsilon : \epsilon \in \{-1,1\}^{\{1,2,3\}}\}$ and has the halving order, see Definition 4.5 and Theorem 4.6. It suffices to show that $\alpha\beta\alpha^{-1}\beta^{-1}$ leaves each of the eight subgroups $U_{\mathfrak{a}}$, $\mathfrak{a} \in \Phi$, invariant. The proof is essentially the same in each case and so only one is described.

Let $\mathfrak{a} = \{\alpha_1, \alpha_2^{-1}, \alpha_3\}$. Then $\alpha\beta\alpha^{-1}\beta^{-1}$ is in $\mathfrak{a}$ and $\alpha\beta\alpha^{-1}\beta^{-1}(U_{\mathfrak{a}}) \geq U_{\mathfrak{a}}$. We compute the index $[\alpha\beta\alpha^{-1}\beta^{-1}(U_{\mathfrak{a}}) : U_{\mathfrak{a}}] = s_{\mathfrak{a}}(\alpha\beta\alpha^{-1}\beta^{-1})$ in terms of the scales of $\alpha$ and $\beta^{-1}$ relative to $U_{\mathfrak{a}}$.

First, since $\beta^{-1}$ is in $\mathfrak{a}$, we have $\beta^{-1}(U_\mathfrak{a}) \geq U_\mathfrak{a}$ and

(28) $$[\beta^{-1}(U_\mathfrak{a}) : U_\mathfrak{a}] = s_\mathfrak{a}(\beta^{-1}).$$

Second, $\beta^{-1}(U)$ is tidy for $\mathfrak{a}$ and Corollary 4.9 shows that $\beta^{-1}(U_\mathfrak{a}) = \beta^{-1}(U)_\mathfrak{a}$. Hence, $\alpha$ being in $\mathfrak{a}$, $\alpha^{-1}\beta^{-1}(U_\mathfrak{a}) \leq \beta^{-1}(U_\mathfrak{a})$ and we have

(29) $$[\beta^{-1}(U_\mathfrak{a}) : \alpha^{-1}\beta^{-1}(U_\mathfrak{a})] = s_\mathfrak{a}(\alpha).$$

(Note that $[\beta^{-1}(U_\mathfrak{a}) : \alpha^{-1}\beta^{-1}(U_\mathfrak{a})] = [\alpha\beta^{-1}(U_\mathfrak{a}) : \beta^{-1}(U_\mathfrak{a})]$ because $\alpha$ is an automorphism.) Third, $\alpha^{-1}\beta^{-1}(U)$ is tidy for $\mathfrak{a}$ and $\alpha^{-1}\beta^{-1}(U_\mathfrak{a}) = \alpha^{-1}\beta^{-1}(U)_\mathfrak{a}$. Hence, $\beta^{-1}$ being in $\mathfrak{a}$, $\beta\alpha^{-1}\beta^{-1}(U_\mathfrak{a}) \leq \alpha^{-1}\beta^{-1}(U_\mathfrak{a})$ and we have

(30) $$[\alpha^{-1}\beta^{-1}(U_\mathfrak{a}) : \beta\alpha^{-1}\beta^{-1}(U_\mathfrak{a})] = s_\mathfrak{a}(\beta^{-1}).$$

Finally, similar reasoning yields that $\alpha\beta\alpha^{-1}\beta^{-1}(U_\mathfrak{a}) \geq \beta\alpha^{-1}\beta^{-1}(U_\mathfrak{a})$ and

(31) $$[\alpha\beta\alpha^{-1}\beta^{-1}(U_\mathfrak{a}) : \beta\alpha^{-1}\beta^{-1}(U_\mathfrak{a})] = s_\mathfrak{a}(\alpha).$$

Since the groups $U_\mathfrak{a}$, $\beta^{-1}(U_\mathfrak{a})$, $\alpha^{-1}\beta^{-1}(U_\mathfrak{a})$, $\beta\alpha^{-1}\beta^{-1}(U_\mathfrak{a})$ and $\alpha\beta\alpha^{-1}\beta^{-1}(U_\mathfrak{a})$ are commensurable, equations (28)–(31) imply that

$$s_\mathfrak{a}(\alpha\beta\alpha^{-1}\beta^{-1}) = \frac{[\alpha\beta\alpha^{-1}\beta^{-1}(U_\mathfrak{a}) : \beta\alpha^{-1}\beta^{-1}(U_\mathfrak{a})][\beta^{-1}(U_\mathfrak{a}) : U_\mathfrak{a}]}{[\beta^{-1}(U_\mathfrak{a}) : \alpha^{-1}\beta^{-1}(U_\mathfrak{a})][\alpha^{-1}\beta^{-1}(U_\mathfrak{a}) : \beta\alpha^{-1}\beta^{-1}(U_\mathfrak{a})]} = 1.$$

Therefore $\alpha\beta\alpha^{-1}\beta^{-1}(U_\mathfrak{a}) = U_\mathfrak{a}$.

Each of the other seven subgroups in (27) may be shown to be invariant under $\alpha\beta\alpha^{-1}\beta^{-1}$ in exactly the same way. Since $U$ is the product of these invariant subgroups, $U$ is invariant under $\alpha\beta\alpha^{-1}\beta^{-1}$. □

**Theorem 4.15.** *Let $\mathfrak{H}$ be a subgroup of $\mathrm{Aut}(G)$ and suppose that $U$ is tidy for $\mathfrak{H}$. Let $\mathfrak{N}_U = \{\alpha \in \mathfrak{H} \mid \alpha(U) = U\}$. Then $\mathfrak{N}_U$ is normal in $\mathfrak{H}$ and $\mathfrak{H}/\mathfrak{N}_U$ is an abelian, torsion-free group no element of which is infinitely divisible.*

**Proof.** Property **S1** of the scale function shows that

(32) $$\mathfrak{N}_U = \{\alpha \in \mathfrak{H} \mid s(\alpha) = 1 = s(\alpha^{-1})\}.$$

Let $\alpha$ be in $\mathfrak{N}_U$ and $\beta$ in $\mathfrak{H}$. Then $\beta^{-1}(U)$ is tidy for $\mathfrak{H}$ by Lemma 4.3 and so $\alpha\beta^{-1}(U) = \beta^{-1}(U)$. Hence $\beta\alpha\beta^{-1}(U) = U$ and $\mathfrak{N}_U$ is normal in $\mathfrak{H}$. The commutator subgroup of $\mathfrak{H}$ is contained in $\mathfrak{N}_U$ by Theorem 4.14. Hence $\mathfrak{H}/\mathfrak{N}_U$ is an abelian group.

If $\alpha^n$ is in $\mathfrak{N}_U$ for some $n \geq 0$, then, by equation (32), $s(\alpha^n) = 1$. Hence, by property **S2** of the scale function, $s(\alpha) = 1$. It may be seen in the same way that $s(\alpha^{-1}) = 1$. Hence $\alpha$ is in $\mathfrak{N}_U$ and $\mathfrak{H}/\mathfrak{N}_U$ is torsion-free.

Now suppose that for some $\alpha\mathfrak{N}_U$ in $\mathfrak{H}/\mathfrak{N}_U$ there are elements $\beta_k\mathfrak{N}_U$ and positive integers $n_k$ with $n_k \to \infty$ as $k \to \infty$ such that $\alpha\mathfrak{N}_U = \beta_k^{n_k}\mathfrak{N}_U$ for each $k$. Then $s(\alpha) = s(\beta_k^{n_k}) = s(\beta_k)^{n_k}$ for each $k$. Since $s(\alpha)$ and $s(\beta_k)$ are positive integers, it follows that $s(\alpha) = 1$. That $s(\alpha^{-1}) = 1$ follows in the same way and so $\alpha$ is in $\mathfrak{N}_U$. Therefore $\mathfrak{H}/\mathfrak{N}_U$ contains no (non-trivial) infinitely divisible elements. □

## 5. Finitely Generated Groups of Automorphisms

For a group of automorphisms, $\mathfrak{H}$, the property of having local tidy subgroups is formally weaker than that of having tidy subgroups. However, it is shown in Theorem 5.5 below that these properties are equivalent when $\mathfrak{H}$ is finitely generated. Since, as seen in §3, abelian groups have local tidy subgroups, it follows that each finitely generated abelian group of automorphisms has a tidy subgroup.

The idea of the proof of this equivalence is to find a finite generating set, $\mathfrak{f}$, for $\mathfrak{H}$ such that any $U$ invariably tidy for $\mathfrak{f}$ is tidy for $\mathfrak{H}$. Example 3.5 shows that not every generating set will do. To find such a set $\mathfrak{f}$, we begin with generators $\alpha_1$, ..., $\alpha_n$ for $\mathfrak{H}$, a compact open subgroup $U$ which is invariably tidy for $\{\alpha_1, \ldots, \alpha_n\}$ and the corresponding factoring $U = \prod_{\mathfrak{a} \in \Phi} U_\mathfrak{a}$, see Theorem 4.6. If the sequence of automorphisms is extended, to $\alpha_1, \ldots, \alpha_n, \alpha_{n+1}$ say, and $U$ is invariably tidy for the extended set, then the corresponding factoring of $U$ is found by factoring each of the $U_\mathfrak{a}$'s as $U_\mathfrak{a} = U_{\mathfrak{a},\alpha_{n+1}+}U_{\mathfrak{a},\alpha_{n+1}-}$. Lemmas 5.1-5.3 show that the relative scale functions are reduced at each factoring. Since these functions take integer values, it follows that this progressive refinement of the factoring of $U$ terminates eventually. This will complete the first stage in the proof of Theorem 5.5.

**Lemma 5.1.** *Let the compact open group $U$ be invariably tidy for the finite set of automorphisms $\mathfrak{a} \cup \{\beta\}$. Then for each $\alpha \in \mathfrak{a}$*

$$\alpha(U_\mathfrak{a}) = \alpha(U_{\mathfrak{a},\beta+})\left(\alpha(U_\mathfrak{a}) \cap U_{\mathfrak{a},\beta--}\right)$$

*and*

$$[\alpha(U_\mathfrak{a}) \cap U_{\mathfrak{a},\beta--} : U_{\mathfrak{a},\beta-}] = \frac{[\alpha(U_{\mathfrak{a},\beta-}) : U_{\mathfrak{a},\beta-}]}{[\alpha(U_{\mathfrak{a},\beta 0}) : U_{\mathfrak{a},\beta 0}]}. \tag{33}$$

**Proof.** Let $w$ be in $\alpha(U_\mathfrak{a})$. Lemma 4.4(ii) shows that $w = w_+ w_-$, where $w_\pm$ is in $\alpha(U_{\mathfrak{a},\beta\pm}) = \alpha(U_\mathfrak{a})_{\beta\pm}$. We modify $w_-$ so that it belongs to $\alpha(U_\mathfrak{a}) \cap U_{\mathfrak{a},\beta--}$.

The set $\{\beta^k(w_-)\}_{k \geq 0}$ is contained in $\alpha(U_\mathfrak{a})$ and so has an accumulation point, $w_0$ say, in $\alpha(U_\mathfrak{a})_{\beta 0}$. It follows that $\beta^l(w_0)$ is in $\alpha(U_\mathfrak{a})_{\beta 0}$ for all $l$. Since, by Lemma 4.10(i), $U_\mathfrak{a}$ is an open subgroup of $\alpha(U_\mathfrak{a})$ there is an $n$ such that $\beta^n(w_-)$ belongs to $w_0 U_\mathfrak{a}$. Put $w'_+ = w_+ \beta^{-n}(w_0)$ and $w'_- = \beta^{-n}(w_0)^{-1} w_-$. Then $w = w'_+ w'_-$ and $w'_+$ belongs to $\alpha(U_{\mathfrak{a},\beta+})$ because $w_+$ and $\beta^{-n}(w_0)$ do and, similarly, $w'_-$ belongs to $\alpha(U_{\mathfrak{a},\beta-})$. We also have that

$$\beta^n(w'_-) = w_0^{-1}\beta^n(w_-),$$

where the latter belongs to $U_\mathfrak{a} \cap \alpha(U_{\mathfrak{a},\beta-})$ by the choice of $n$. Now Lemma 4.10(i) implies that $U_\mathfrak{a} \cap \alpha(U_{\mathfrak{a},\beta-}) = U_{\mathfrak{a},\beta-}$ and so $w'_-$ belongs to $U_{\mathfrak{a},\beta--}$ as required.

Choosing $w$ in $\alpha(U_{\mathfrak{a},\beta-})$ in the preceding argument shows that

$$\alpha(U_{\mathfrak{a},\beta-}) = \alpha(U_{\mathfrak{a},\beta 0})\left(\alpha(U_\mathfrak{a}) \cap U_{\mathfrak{a},\beta--}\right).$$

Hence

$$[\alpha(U_{\mathfrak{a},\beta-}) : U_{\mathfrak{a},\beta-}] = [\alpha(U_{\mathfrak{a},\beta-}) : \alpha(U_\mathfrak{a}) \cap U_{\mathfrak{a},\beta--}][\alpha(U_\mathfrak{a}) \cap U_{\mathfrak{a},\beta--} : U_{\mathfrak{a},\beta-}] \tag{34}$$

and the map

$$\Theta : \alpha(U_{\mathfrak{a},\beta 0})/U_{\mathfrak{a},\beta 0} \to \alpha(U_{\mathfrak{a},\beta-})/\left(\alpha(U_\mathfrak{a}) \cap U_{\mathfrak{a},\beta--}\right)$$

defined by

$$\Theta(\alpha(u)U_{\mathfrak{a},\beta 0}) = \alpha(u)\left(\alpha(U_\mathfrak{a}) \cap U_{\mathfrak{a},\beta--}\right)$$

is surjective. The map $\Theta$ is well-defined because $U_{\mathfrak{a},\beta 0} \leq \alpha(U_\mathfrak{a}) \cap U_{\mathfrak{a},\beta--}$ and is injective because $\alpha(U_{\mathfrak{a},\beta 0}) \cap U_{\mathfrak{a},\beta--} = U_{\mathfrak{a},\beta 0}$. Therefore

$$[\alpha(U_{\mathfrak{a},\beta-}) : \alpha(U_\mathfrak{a}) \cap U_{\mathfrak{a},\beta--}] = [\alpha(U_{\mathfrak{a},\beta 0}) : U_{\mathfrak{a},\beta 0}]$$

and substitution into (34) yields (33). $\square$

**Lemma 5.2.** *Let the compact open group $U$ be invariably tidy for the finite set of automorphisms $\mathfrak{a} \cup \{\beta\}$. Then for each $\alpha \in \mathfrak{a}$*

$$(35) \qquad [\alpha(U_\mathfrak{a}) : U_\mathfrak{a}] = \frac{[\alpha(U_{\mathfrak{a},\beta+}) : U_{\mathfrak{a},\beta+}][\alpha(U_{\mathfrak{a},\beta-}) : U_{\mathfrak{a},\beta-}]}{[\alpha(U_{\mathfrak{a},\beta 0}) : U_{\mathfrak{a},\beta 0}]}.$$

**Proof.** The group $\alpha(U_\mathfrak{a}) \cap U_{\mathfrak{a},\beta--}$ is compact because $\alpha(U_\mathfrak{a})$ is compact and $U_{\mathfrak{a},\beta--}$ is closed, by Lemma 4.10(iv). Hence, by definition of $U_{\mathfrak{a},\beta--}$, there is an $n$ such that $\beta^n(\alpha(U_\mathfrak{a}) \cap U_{\mathfrak{a},\beta--}) < U_\mathfrak{a}$. Now $\beta^n(\alpha(U_\mathfrak{a}) \cap U_{\mathfrak{a},\beta--}) = \beta^n \alpha(U_\mathfrak{a}) \cap U_{\mathfrak{a},\beta--}$ and Lemma 2.6(i) implies that $\beta^n \alpha(U_\mathfrak{a}) \cap U_{\mathfrak{a},\beta--} = \beta^n \alpha(U_{\mathfrak{a},\beta-})$. Hence

$$(36) \qquad \beta^n(\alpha(U_\mathfrak{a}) \cap U_{\mathfrak{a},\beta--}) = \beta^n \alpha(U_{\mathfrak{a},\beta-}) < U_\mathfrak{a}.$$

Put $U_{\mathfrak{a},n} = \bigcap_{k=0}^n \beta^k(U_\mathfrak{a})$. Then $U_{\mathfrak{a},n} = \left(\bigcap_{k=0}^n \beta^k(U)\right)_\mathfrak{a}$, where $\bigcap_{k=0}^n \beta^k(U)$ is invariably tidy for $\mathfrak{a} \cup \{\beta\}$, and so Theorem 4.12 shows that

$$(37) \qquad [\alpha(U_\mathfrak{a}) : U_\mathfrak{a}] = [\alpha(U_{\mathfrak{a},n}) : U_{\mathfrak{a},n}].$$

Furthermore, $\bigcap_{k=0}^n \beta^k(V) = V_{\beta+} \beta^n(V_{\beta-})$ for every $V$ tidy for $\beta$ and so

$$(38) \qquad U_{\mathfrak{a},n} = U_{\mathfrak{a},\beta+} \beta^n(U_{\mathfrak{a},\beta-}).$$

Now, since $U_{\mathfrak{a},n} < \alpha(U_{\mathfrak{a},n}) \cap U_\mathfrak{a} < \alpha(U_{\mathfrak{a},n})$, we have

$$[\alpha(U_{\mathfrak{a},n}) : U_{\mathfrak{a},n}] = [\alpha(U_{\mathfrak{a},n}) : \alpha(U_{\mathfrak{a},n}) \cap U_\mathfrak{a}][\alpha(U_{\mathfrak{a},n}) \cap U_\mathfrak{a} : U_{\mathfrak{a},n}]$$

and so (35) will follow from (37) once it is established that

$$(39) \qquad [\alpha(U_{\mathfrak{a},n}) : \alpha(U_{\mathfrak{a},n}) \cap U_\mathfrak{a}] = [\alpha(U_{\mathfrak{a},\beta+}) : U_{\mathfrak{a},\beta+}]$$

and

$$(40) \qquad [\alpha(U_{\mathfrak{a},n}) \cap U_\mathfrak{a} : U_{\mathfrak{a},n}] = \frac{[\alpha(U_{\mathfrak{a},\beta-}) : U_{\mathfrak{a},\beta-}]}{[\alpha(U_{\mathfrak{a},\beta 0}) : U_{\mathfrak{a},\beta 0}]}.$$

To establish (39) we show that the map

$$\Theta : \alpha(U_{\mathfrak{a},\beta+})/U_{\mathfrak{a},\beta+} \to \alpha(U_{\mathfrak{a},n})/(\alpha(U_{\mathfrak{a},n}) \cap U_\mathfrak{a})$$

defined by $\Theta(\alpha(u)U_{\mathfrak{a},\beta+}) = \alpha(u)(\alpha(U_{\mathfrak{a},n}) \cap U_\mathfrak{a})$ is a bijection. This map is well-defined because $U_{\mathfrak{a},\beta+}$ is a subgroup of $\alpha(U_{\mathfrak{a},n}) \cap U_\mathfrak{a}$. By (38),

$$\alpha(U_{\mathfrak{a},n}) = \alpha(U_{\mathfrak{a},\beta+})\alpha\beta^n(U_{\mathfrak{a},\beta-})$$

and the choice of $n$, see (36), guarantees that

$$\alpha(U_{\mathfrak{a},n}) \cap U_\mathfrak{a} = U_{\mathfrak{a},\beta+}\alpha\beta^n(U_{\mathfrak{a},\beta-}).$$

Hence the map $\Theta$ is surjective. It is injective because

$$\alpha(U_{\mathfrak{a},\beta+}) \cap (\alpha(U_{\mathfrak{a},n}) \cap U_\mathfrak{a}) = U_{\mathfrak{a},\beta+}.$$

To establish (40) it suffices, by Lemma 5.1, to show that $[\alpha(U_{\mathfrak{a},n}) \cap U_\mathfrak{a} : U_{\mathfrak{a},n}] = [\alpha(U_\mathfrak{a}) \cap U_{\mathfrak{a},\beta--} : U_{\mathfrak{a},\beta-}]$ and for this it suffices, since $\beta$ is an automorphism, to show that $[\alpha(U_{\mathfrak{a},n}) \cap U_{\mathfrak{a}} : U_{\mathfrak{a},n}] = [\beta^n (\alpha(U_{\mathfrak{a}}) \cap U_{\mathfrak{a},\beta--}) : \beta^n (U_{\mathfrak{a},\beta-})]$. To this end, define a map

$$\Upsilon : \beta^n (\alpha(U_{\mathfrak{a}}) \cap U_{\mathfrak{a},\beta--}) / \beta^n (U_{\mathfrak{a},\beta-}) \to (\alpha(U_{\mathfrak{a},n}) \cap U_{\mathfrak{a}}) / U_{\mathfrak{a},n}$$

by $\Upsilon(w\beta^n (U_{\mathfrak{a},\beta-})) = wU_{\mathfrak{a},n}$. Then $\Upsilon$ is well-defined because $\beta^n (U_{\mathfrak{a},\beta-}) < U_{\mathfrak{a},n}$ and because $\beta^n (\alpha(U_{\mathfrak{a}}) \cap U_{\mathfrak{a},\beta--}) < \alpha(U_{\mathfrak{a},n}) \cap U_{\mathfrak{a}}$ by the choice of $n$, see (36). Now

$$U_{\mathfrak{a},n} \cap \beta^n (\alpha(U_{\mathfrak{a}}) \cap U_{\mathfrak{a},\beta--}) = U_{\mathfrak{a},\beta+}\beta^n(U_{\mathfrak{a},\beta-}) \cap \alpha\beta^n(U_{\mathfrak{a},\beta-}) = \beta^n (U_{\mathfrak{a},\beta-})$$

and so $\Upsilon$ is injective. Finally, since $\alpha(U_{\mathfrak{a},n}) \cap U_{\mathfrak{a}} = \alpha\beta^n(U_{\mathfrak{a},\beta-})U_{\mathfrak{a},\beta+}$, (36) shows that $\Upsilon$ is surjective. □

**Lemma 5.3.** *Let $\mathfrak{a} \subseteq \mathfrak{b}$ be finite sets of automorphisms of $G$, let $\beta \in \langle \mathfrak{a} \rangle$ and suppose that $\beta = \eta \prod_{\alpha \in \mathfrak{a}} \alpha^{r_\alpha}$, where $\eta$ is in the commutator subgroup of $\langle \mathfrak{a} \rangle$. Let $U$ be a compact open subgroup of $G$ which is invariably tidy for $\mathfrak{b} \cup \{\beta, \eta\}$ and suppose that $\beta(U_{\mathfrak{b}}) \not\geq U_{\mathfrak{b}}$. Then there is $\alpha \in \mathfrak{a}$ such that*

$$[\alpha(U_{\mathfrak{b},\beta+}) : U_{\mathfrak{b},\beta+}] < [\alpha(U_{\mathfrak{b}}) : U_{\mathfrak{b}}].$$

**Proof.** Lemma 5.2 shows that, if the claim fails, then

(41) $\qquad [\alpha(U_{\mathfrak{b},\beta-}) : U_{\mathfrak{b},\beta-}] = [\alpha(U_{\mathfrak{b},\beta 0}) : U_{\mathfrak{b},\beta 0}]$ for each $\alpha \in \mathfrak{a}$.

Since $\eta$ is in the commutator subgroup of $\langle \mathfrak{a} \rangle$, it follows from Theorem 4.14 that $\eta(V) = V$ for any $V$ tidy for $\eta$. Hence, by the invariable tidiness of $U$, we have $\beta(U) = \left( \prod_{\alpha \in \mathfrak{a}} \alpha^{r_\alpha} \right) (U)$ and

$$[U_{\mathfrak{b},\beta-} : \beta(U_{\mathfrak{b},\beta-})] = [U_{\mathfrak{b},\beta-} : \left( \prod_{\alpha \in \mathfrak{a}} \alpha^{r_\alpha} \right) (U_{\mathfrak{b},\beta-})].$$

An argument like that used in the proof of Theorem 4.14 shows that

$$[U_{\mathfrak{b},\beta-} : \left( \prod_{\alpha \in \mathfrak{a}} \alpha^{r_\alpha} \right) (U_{\mathfrak{b},\beta-})] = \prod_{\alpha \in \mathfrak{a}} [\alpha(U_{\mathfrak{b},\beta-}) : U_{\mathfrak{b},\beta-}]^{-r_\alpha}$$

and, by (41), the right hand side equals

$$\prod_{\alpha \in \mathfrak{a}} [\alpha(U_{\mathfrak{b},\beta 0}) : U_{\mathfrak{b},\beta 0}]^{-r_\alpha}.$$

The same argument as above shows that this expression equals $[U_{\mathfrak{b},\beta 0} : \beta(U_{\mathfrak{b},\beta 0})]$, which equals 1. Hence $U_{\mathfrak{b},\beta-} = U_{\mathfrak{b},\beta 0}$ and $\beta(U_{\mathfrak{b}}) \geq U_{\mathfrak{b}}$, in contradiction to the hypothesis. □

**Lemma 5.4.** *Let $\mathfrak{H}$ be a finitely generated group of automorphisms of $G$ and let $\mathfrak{L}$ be a normal subgroup of $\mathfrak{H}$ such that $\mathfrak{H}/\mathfrak{L}$ is abelian. Then there are elements $\gamma_1$, …, $\gamma_r$ in $\mathfrak{L}$ such that*

$$\mathfrak{L} = \left\langle \alpha \gamma_j \alpha^{-1} : \alpha \in G, j \in \{1, \ldots, r\} \right\rangle.$$

**Proof.** Let $\mathfrak{H}'$ be the commutator subgroup of $\mathfrak{H}$. Then $\mathfrak{H}' \leq \mathfrak{L}$ and $\mathfrak{L}/\mathfrak{H}'$ is a finitely generated abelian group. Suppose that $\alpha_1$, …, $\alpha_n$ generate $\mathfrak{H}$. Then every commutator in $\mathfrak{H}$ is a product of conjugates of commutators of the form $[\alpha_k, \alpha_l]$, where $\alpha_k$ and $\alpha_l$ are in the generating set. Take the $\gamma_j$'s to be these commutators. □

**Theorem 5.5.** *Let $\mathfrak{H}$ be be a finitely generated subgroup of $\operatorname{Aut}(G)$ which has local tidy subgroups. Then there is a compact open subgroup $U$ which is tidy for $\mathfrak{H}$.*

**Proof.** Suppose that $\alpha_1, \ldots, \alpha_n$ are automorphisms with $\mathfrak{H} = \langle \alpha_1, \ldots, \alpha_n \rangle$. Let $U$ be a compact open subgroup which is invariably tidy for $\{\alpha_1, \ldots, \alpha_n\}$ and let

$$U = \prod_{\mathfrak{a} \in \Phi} U_{\mathfrak{a}} \tag{42}$$

be the factoring of $U$ corresponding to the sequence $\alpha_1, \ldots, \alpha_n$.

If the sequence of automorphisms is extended, to $\alpha_1, \ldots, \alpha_n, \alpha_{n+1}, \ldots, \alpha_{n+p}$ say, and $U$ is invariably tidy for $\{\alpha_1, \ldots, \alpha_n, \alpha_{n+1}, \ldots, \alpha_{n+p}\}$, then the factoring of $U$ corresponding to the extended sequence is a refinement of (42) as follows. For each function $\epsilon$ in $\{-1, 1\}^{\{1,\ldots,p\}}$ put $\mathfrak{b}_\epsilon = \left\{ \alpha_{n+j}^{\epsilon(j)} : j \in \{1, \ldots, p\} \right\}$ and define

$$\Psi = \left\{ \mathfrak{b}_\epsilon : \epsilon \in \{-1, 1\}^{\{1,\ldots,p\}} \right\}.$$

Then for each $\mathfrak{a}$ in $\Phi$ we have, by Lemma 4.4, $U_{\mathfrak{a}} = \prod_{\mathfrak{b} \in \Psi} U_{\mathfrak{a} \cup \mathfrak{b}}$ and

$$U = \prod_{\mathfrak{a} \in \Phi} U_{\mathfrak{a}} = \prod_{\mathfrak{a} \in \Phi} \left( \prod_{\mathfrak{b} \in \Psi} U_{\mathfrak{a} \cup \mathfrak{b}} \right) = \prod_{(\mathfrak{a}, \mathfrak{b}) \in \Phi \times \Psi} U_{\mathfrak{a} \cup \mathfrak{b}} = \prod_{\mathfrak{d} \in \Phi'} U_{\mathfrak{d}} \tag{43}$$

is the factoring of $U$ corresponding to the sequence of automorphisms $\alpha_1, \ldots, \alpha_n$, $\alpha_{n+1}, \ldots, \alpha_{n+p}$. The first stage of the proof extends the sequence of automorphisms to the point where every possible choice for $\alpha_{n+p+1}$ satisfies, for each factor $U_{\mathfrak{a} \cup \mathfrak{b}}$, either $\alpha_{n+p+1}(U_{\mathfrak{a} \cup \mathfrak{b}}) \geq U_{\mathfrak{a} \cup \mathfrak{b}}$ or $\alpha_{n+p+1}(U_{\mathfrak{a} \cup \mathfrak{b}}) \leq U_{\mathfrak{a} \cup \mathfrak{b}}$.

Suppose that is not already so for the sequence $\alpha_1, \ldots, \alpha_n$. Then there is $\beta$ in $\mathfrak{H}$ and $\mathfrak{a}$ in $\Phi$, and there is $U$ invariably tidy for $\{\alpha_1, \ldots, \alpha_n, \beta\}$, such that $\beta(U_{\mathfrak{a}}) \not\geq U_{\mathfrak{a}}$ and $\beta(U_{\mathfrak{a}}) \not\leq U_{\mathfrak{a}}$, whence $U_{\mathfrak{a},\beta+}$ and $U_{\mathfrak{a},\beta-}$ are proper subgroups of $U_{\mathfrak{a}}$. Lemma 5.3 shows that, for any $\mathfrak{c}$ in $\Phi$ with $U_{\mathfrak{c},\beta+}$ a proper subgroup of $U_{\mathfrak{c}}$, we have $s_{\mathfrak{c} \cup \{\beta\}}(\alpha) < s_{\mathfrak{c}}(\alpha)$ for some $\alpha$ in $\mathfrak{c}$ and similarly if $U_{\mathfrak{c},\beta-}$ is a proper subgroup of $U_{\mathfrak{c}}$. (Note that, since relative scales are independent of the tidy subgroup, the Lemma applies even though the present subgroup $U$ might not be tidy for $\eta$.) Choose $\alpha_{n+1} = \beta$. Then, taking $p = 1$, the refinement (43) satisfies that $s_{\mathfrak{c} \cup \mathfrak{b}}(\alpha) < s_{\mathfrak{c}}(\alpha)$ for some $\alpha$ in $\mathfrak{c}$ whenever $U_{\mathfrak{c} \cup \mathfrak{b}} \neq U_{\mathfrak{c}}$.

Continue, choosing $\alpha_{n+2}$ *etc.*, in the same way so that, whenever a new subgroup appears as a factor in (43) the relative scale of some $\alpha$ in $\{\alpha_j, \alpha_j^{-1} : j \in \{1, \ldots, n\}\}$ is strictly reduced on that factor. These relative scales are being reduced from the initial values of $s_{\mathfrak{a}}(\alpha)$, where $\mathfrak{a} \in \Phi$ and $\alpha \in \mathfrak{a}$. Since $\{s_{\mathfrak{a}}(\alpha) : \alpha \in \mathfrak{a}, \mathfrak{a} \in \Phi\}$ is a finite set of positive integers, a point is reached, with $\alpha_{n+p}$ say, where no further reduction in relative scales is possible. At that point we have for every $\beta$ in $\mathfrak{H}$ that, if $U$ is invariably tidy for $\{\alpha_1, \ldots, \alpha_{n+p}, \beta\}$, then for every factor $U_{\mathfrak{a} \cup \mathfrak{b}}$ in (43) either $\beta(U_{\mathfrak{a} \cup \mathfrak{b}}) \geq U_{\mathfrak{a} \cup \mathfrak{b}}$ or $\beta(U_{\mathfrak{a} \cup \mathfrak{b}}) \leq U_{\mathfrak{a} \cup \mathfrak{b}}$. Now fix the sequence of automorphisms $\alpha_1, \ldots, \alpha_n, \alpha_{n+1}, \ldots, \alpha_{n+p}$ and let $\Phi'$ be the index set for the corresponding factoring of any $U$ invariably tidy for $\{\alpha_1, \ldots, \alpha_{n+p}\}$, see Theorem 4.6.

The next stage of the proof is to show that there is a $U$ such that, for every $\mathfrak{d}$ in $\Phi'$ and $\beta$ in $\mathfrak{H}$, either $\beta(U_{\mathfrak{d}}) \geq U_{\mathfrak{d}}$ or $\beta(U_{\mathfrak{d}}) \leq U_{\mathfrak{d}}$. For each $\mathfrak{d}$ in $\Phi'$ define a map $\psi_{\mathfrak{d}} : \mathfrak{H} \to \mathbf{Q}^+$ by

$$\psi_{\mathfrak{d}}(\beta) = \begin{cases} [\beta(U_{\mathfrak{d}}) : U_{\mathfrak{d}}], & \text{if } \beta(U_{\mathfrak{d}}) \geq U_{\mathfrak{d}} \\ [U_{\mathfrak{d}} : \beta(U_{\mathfrak{d}})]^{-1}, & \text{if } \beta(U_{\mathfrak{d}}) \leq U_{\mathfrak{d}}, \end{cases} \tag{44}$$

where $U$ is a subgroup invariably tidy for $\{\alpha_1, \ldots, \alpha_{n+p}, \beta\}$. Theorem 4.12 shows that $\psi_{\mathfrak{d}}$ is well-defined. For $\beta$ and $\gamma$ in $\mathfrak{H}$ and expanding on $U_{\mathfrak{d}}$ we have

$$[\beta\gamma(U_{\mathfrak{d}}) : U_{\mathfrak{d}}] = [\beta(\gamma(U)_{\mathfrak{d}}) : \gamma(U)_{\mathfrak{d}}][\gamma(U_{\mathfrak{d}}) : U_{\mathfrak{d}}],$$

and so, since $\gamma(U)$ is tidy for $\beta$ if $U$ is, $\psi_{\mathfrak{d}}(\beta\gamma) = \psi_{\mathfrak{d}}(\beta)\psi_{\mathfrak{d}}(\gamma)$. By checking other cases in the same way it may be seen that each $\psi_{\mathfrak{d}}$ is a homomorphism from $\mathfrak{H}$ to the multiplicative group of positive rationals. If $\beta$ and $\gamma$ are in $\mathfrak{H}$ and $\psi_{\mathfrak{d}}(\beta\gamma^{-1}) \geq 1$, then $\psi_{\mathfrak{d}}(\beta)\psi_{\mathfrak{d}}(\gamma)^{-1} = \psi_{\mathfrak{d}}(\beta\gamma^{-1})$ is a positive integer. It follows that for each $\mathfrak{d}$ in $\Phi'$ there is a positive integer, $t_{\mathfrak{d}}$, such that $\psi_{\mathfrak{d}}(\mathfrak{H}) = \{t_{\mathfrak{d}}^m : m \in \mathbf{Z}\}$.

For each $\mathfrak{d}$ in $\Phi'$, choose an automorphism $\gamma_{\mathfrak{d}}$ in $\mathfrak{H}$ such that $\psi_{\mathfrak{d}}(\gamma_{\mathfrak{d}}) = t_{\mathfrak{d}}$. Then

$$\mathfrak{H} = \{\gamma_{\mathfrak{d}}^m : m \in \mathbf{Z}\} \ltimes \ker\psi_{\mathfrak{d}},$$

is the semidirect product of an infinite cyclic group and the group of automorphisms which leave $U_{\mathfrak{d}}$ invariant. Choose automorphisms $\eta_{\mathfrak{d}}^{(1)}, \ldots, \eta_{\mathfrak{d}}^{(r_{\mathfrak{d}})}$ such that

$$\ker\psi_{\mathfrak{d}} = \left\langle \alpha\eta_{\mathfrak{d}}^{(j)}\alpha^{-1} : \alpha \in \mathfrak{H}, j \in \{1, \ldots, r_{\mathfrak{d}}\}\right\rangle,$$

see Lemma 5.4. Now let $U$ be invariably tidy for

(45) $\quad\quad\quad \{\alpha_1, \ldots, \alpha_p\} \cup \{\gamma_{\mathfrak{d}} : \mathfrak{d} \in \Phi'\} \cup \{\eta_{\mathfrak{d}}^{(j)} : \mathfrak{d} \in \Phi', 1 \leq j \leq r_{\mathfrak{d}}\}.$

Then $\eta_{\mathfrak{d}}^{(j)}(U_{\mathfrak{d}}) = U_{\mathfrak{d}}$ for each $\mathfrak{d}$ in $\Phi'$ and $j \in \{1, \ldots, r_{\mathfrak{d}}\}$. Since $U$ is invariably tidy, it follows that $\alpha\eta_{\mathfrak{d}}^{(j)}\alpha^{-1}(U_{\mathfrak{d}}) = U_{\mathfrak{d}}$ for each $j \in \{1, \ldots, r_{\mathfrak{d}}\}$ and $\alpha$ in $\mathfrak{H}$. Since the elements $\alpha\eta_{\mathfrak{d}}^{(j)}\alpha^{-1}$, $j \in \{1, \ldots, r_{\mathfrak{d}}\}$, $\alpha \in \mathfrak{H}$, generate $\ker\psi_{\mathfrak{d}}$, it follows that $\eta(U_{\mathfrak{d}}) = U_{\mathfrak{d}}$ for each $\eta$ in $\ker\psi_{\mathfrak{d}}$. Now let $\gamma$ be in $\mathfrak{H}$. For each $\mathfrak{d}$ we may write $\gamma = \gamma_{\mathfrak{d}}^m \eta_{\mathfrak{d}}$, where $\eta_{\mathfrak{d}}$ is in $\ker\psi_{\mathfrak{d}}$. Hence

$$\gamma(U_{\mathfrak{d}}) = \gamma_{\mathfrak{d}}^m \eta_{\mathfrak{d}}(U_{\mathfrak{d}}) = \gamma_{\mathfrak{d}}^m(U_{\mathfrak{d}})$$

and it follows, since $U$ is tidy for $\gamma_{\mathfrak{d}}$, that either $\gamma(U_{\mathfrak{d}}) \geq U_{\mathfrak{d}}$ or $\gamma(U_{\mathfrak{d}}) \leq U_{\mathfrak{d}}$.

The final stage of the proof enlarges the finite set in (45) still further so that any subgroup $U$ which is invariably tidy for the enlarged set if tidy for $\mathfrak{H}$. Each $\gamma$ in $\mathfrak{H}$ partitions $\Phi'$ into three subsets

(46) $\quad\quad\quad\quad\quad\quad\quad \Phi'_{\gamma+} := \{\mathfrak{d} \in \Phi' : \gamma(U_{\mathfrak{d}}) > U_{\mathfrak{d}}\},$

(47) $\quad\quad\quad\quad\quad\quad\quad \Phi'_{\gamma-} := \{\mathfrak{d} \in \Phi' : \gamma(U_{\mathfrak{d}}) < U_{\mathfrak{d}}\}$

(48) $\quad\quad\quad\quad \text{and}\quad \Phi'_{\gamma 0} := \{\mathfrak{d} \in \Phi' : \gamma(U_{\mathfrak{d}}) = U_{\mathfrak{d}}\}.$

Now Lemma 4.4, applied with $\beta$ equal to $\alpha_1, \ldots, \alpha_{n+p}$ in turn implies that

$$U_{\gamma+} = \prod_{\mathfrak{d} \in \Phi'} U_{\gamma+,\mathfrak{d}},$$

where $U_{\gamma+,\mathfrak{d}} = U_{\gamma+} \cap U_{\mathfrak{d}}$. (Recall that we do not need $U$ tidy for $\gamma$ for this conclusion.) It is clear that if $\mathfrak{d}$ is in $\Phi'_{\gamma+}$ or $\Phi'_{\gamma 0}$, then $U_{\gamma+,\mathfrak{d}} = U_{\mathfrak{d}}$. On the other hand, if $\mathfrak{d}$ is in $\Phi'_{\gamma-}$, then $\gamma_{\mathfrak{d}}(U_{\gamma+,\mathfrak{d}}) = U_{\gamma+,\mathfrak{d}}$. Since $\mathfrak{H} = \langle\gamma_{\mathfrak{d}}\rangle\ker\psi_{\mathfrak{d}}$, it follows that $\alpha(U_{\gamma+,\mathfrak{d}}) = U_{\gamma+,\mathfrak{d}}$ for every $\alpha$ in $\mathfrak{H}$. Hence $U_{\gamma+,\mathfrak{d}} = \bigcap_{\alpha \in \mathfrak{H}}\alpha(U)$ and is contained in every other factor in the product. Similar arguments apply to $U_{\gamma-}$ and so

(49) $\quad\quad U_{\gamma+} = \prod_{\mathfrak{d} \in \Phi'_{\gamma+} \cup \Phi'_{\gamma 0}} U_{\mathfrak{d}} \quad \text{and} \quad U_{\gamma-} = \prod_{\mathfrak{d} \in \Phi'_{\gamma-} \cup \Phi'_{\gamma 0}} U_{\mathfrak{d}},$

where the order of the products is given by the halving order on $\Phi'$.

Since $\Phi'$ is finite, it can be partitioned as $\Phi' = \Phi'_{\gamma+} \cup \Phi'_{\gamma-} \cup \Phi'_{\gamma 0}$ in only finitely many ways. Let $q$ be the number of such partitions and choose automorphisms $\gamma_1, \ldots, \gamma_q$ which give rise to them. Let $U$ be a compact open subgroup which is invariably tidy for

$$\{\alpha_1, \ldots, \alpha_p\} \cup \{\gamma_{\mathfrak{d}} : \mathfrak{d} \in \Phi'\} \cup \{\eta_{\mathfrak{d}}^{(j)} : \mathfrak{d} \in \Phi', \ 1 \leq j \leq r_{\mathfrak{d}}\} \cup \{\gamma_1, \ldots, \gamma_q\}.$$

Let $\alpha$ be in $\mathfrak{H}$. Then there is $\gamma_j$ such that $\Phi'_{\alpha+} = \Phi'_{\gamma_j+}$, $\Phi'_{\alpha-} = \Phi'_{\gamma_j-}$ and $\Phi'_{\alpha 0} = \Phi'_{\gamma_j 0}$, whence (49) implies that $U_{\alpha+} = U_{\gamma_j+}$, $U_{\alpha-} = U_{\gamma_j-}$ and $U_{\alpha 0} = U_{\gamma_j 0}$. Since $U$ satisfies $\mathbf{T1}(\gamma_j)$, it follows that

$$U = U_{\alpha+} U_{\alpha-},$$

that is, $U$ satisfies $\mathbf{T1}(\alpha)$. Finally, (49) implies that there are integers $k$ and $l$ such that $\gamma_j(U_{\gamma_j+}) \leq \alpha^k(U_{\alpha+}) \leq \gamma_j^l(U_{\gamma_j+})$. Hence

$$\bigcup_{k \geq 0} \alpha^k(U_{\alpha+}) = \bigcup_{l \geq 0} \gamma_j^k(U_{\gamma_j+}),$$

which is closed because $U$ satisfies $\mathbf{T2}(\gamma_j)$. Therefore $U$ satisfies $\mathbf{T2}(\alpha)$ and we have shown that $U$ is tidy for the general element $\alpha$ of $\mathfrak{H}$. $\square$

If $U$ is tidy for $\mathfrak{H}$, then it is invariably tidy for the automorphisms chosen in the proof of the Theorem.

**Corollary 5.6.** *Let $\mathfrak{H}$ be a finitely generated subgroup of $\mathrm{Aut}(G)$ and $U$ be tidy for $\mathfrak{H}$. Then there are automorphisms $\alpha_1, \ldots, \alpha_r$ in $\mathfrak{H}$ such that the corresponding decomposition,*

$$U = \prod_{\mathfrak{a} \in \Phi} U_{\mathfrak{a}},$$

*satisfies that for every $\alpha$ in $\mathfrak{H}$ and $\mathfrak{a}$ in $\Phi$ either $\alpha(U_{\mathfrak{a}}) \geq U_{\mathfrak{a}}$ or $\alpha(U_{\mathfrak{a}}) \leq U_{\mathfrak{a}}$.*

The automorphisms $\alpha_1$ and $\alpha_2$ and the subgroup $U$ in Example 3.5 are already at the first stage of the argument in the proof of Theorem 5.5. The factoring, in this case $U = U_{\alpha_1} = U_{\alpha_2}$, cannot be further refined. However the second stage of the proof does identify a subgroup tidy for $\langle \alpha_1, \alpha_2 \rangle$ in this case. It may be that a subgroup tidy for $\mathfrak{H}$ is always produced at this stage of the argument. I am not aware of an example where the third stage of the proof is required.

The hypothesis in Theorem 5.5 that $\mathfrak{H}$ be finitely generated is necessary.

**Example 5.7.** Let $S_3$ be the permutation group on three symbols. Denote the identity of $S_3$ by $e$, the three order 2 elements by $\sigma_j$, $j = 1, 2, 3$ and the order 3 elements by $\tau$ and $\tau^2$. Define

$$G = \left\{ f \in S_3^{\mathbf{Z}} : f(k) \in \{e, \sigma_1\} \text{ for all but finitely many } k \right\}.$$

Then $G$ is a group under coordinatewise multiplication. For each natural number $n$ define the subgroup

$$G_n = \{f \in G : f(m) \in \{e, \sigma_1\} \text{ for every } m \text{ and } f(m) = e \text{ for } -n < m < n\}.$$

Note that $G_0 = \{f \in S_3^{\mathbf{Z}} : f(k) \in \{e, \sigma_1\} \text{ for every } k\}$. Define a topology on $G$ by letting $\{G_n : n \in \mathbf{N}\}$ be a base of neigbourhoods for $e_G$. Then $G$ is a totally disconnected locally compact group. The $n$th coordinate projection map, $f \mapsto f(n)$, of $G$ onto $S_3$ will be denoted $\pi_n$.

Define, for each $n$ in $\mathbf{Z}$, an automorphism $\alpha_n$ of $G$ by

$$\alpha_n(f)(k) = \begin{cases} \tau f(k)\tau^{-1}, & \text{if } k = n, \\ f(k), & \text{otherwise,} \end{cases} \qquad (f \in G).$$

Then $s(\alpha_n) = 1$ for each $n$ and any subgroup tidy for $\alpha_n$ satisfies $\alpha_n(U) = U$. Since $\{e, \sigma_1\}$ is not normal in $S_3$, a compact open subgroup $U$ is tidy for $\alpha_n$ if and only if $\pi_n(U) = \{e\}$, $\{e, \tau, \tau^2\}$ or $S_3$. Let $\mathfrak{H} = \langle \alpha_n : n \in \mathbf{Z} \rangle$.

Now for each $n$ the subgroup $G_n$ is invariably tidy for $\{\alpha_k : -n < k < n\}$. Since $\mathfrak{H} = \bigcup_{n \in \mathbf{N}} \langle \alpha_k : -n < k < n \rangle$, it follows that $\mathfrak{H}$ has local tidy subgroups. However $\mathfrak{H}$ is not finitely generated and any compact group $U$ invariant for all $\alpha$ in $\mathfrak{H}$ must satisfy $\pi_n(U) = \{e\}$ for all but finitely many $n$. Hence $\mathfrak{H}$ does not have tidy subgroups.

Define another automorphism $\alpha$ of $G$ by

$$\alpha(f)(k) = f(k+1) \qquad (f \in G),$$

and put $\mathfrak{H}_1 = \langle \mathfrak{H}, \alpha \rangle$. Then $\mathfrak{H}_1$ is finitely generated. Now $G_0$ is tidy for $\alpha$, and is the only such subgroup, but it is not tidy for $\alpha_0$. Hence there is no subgroup tidy for $\{\alpha, \alpha_0\}$ and $\mathfrak{H}_1$ does not have local tidy subgroups.

A compact open subgroup tidy for automorphisms $\alpha$ and $\beta$ is not, in general, tidy for their product. That is so in the following special case however.

**Proposition 5.8.** *Let $\alpha$ and $\beta$ be commuting automorphisms of $G$ and let $U$ be a compact open subgroup of $G$ which is tidy for $\alpha$ and satisfies $\beta(U) = U$. Then $U$ is tidy for $\alpha\beta$.*

**Proof.** Since $(\alpha\beta)^n(U) = \alpha^n\beta^n(U) = \alpha^n(U)$, we have $U_{\alpha\beta\pm} = U_{\alpha\pm}$. It is also easy to see that $\beta(U_{\alpha\pm}) = U_{\alpha\pm}$, and so $U_{\alpha\beta\pm\pm} = U_{\alpha\pm\pm}$. □

For an element $x$ of $G$, we say that the compact open subgroup $U$ is *tidy for $x$* if $U$ is tidy for the inner automorphism $\alpha_x : y \mapsto xyx^{-1}$. The following topological version of Theorem 5.5 for abelian groups inner automorphisms may be proved with the aid of Proposition 5.8.

**Theorem 5.9.** *Let $H$ be a compactly generated closed abelian subgroup of $G$. Then there is a subgroup $U$ tidy for $H$.*

**Proof.** Let $V$ be a compact open subgroup of $G$. Then $H \cap V$ is a compact open normal subgroup of $H$. Since $H$ is compactly generated, there are $x_1, \ldots, x_n$ in $H$ such that $H = \langle x_1, \ldots, x_n \rangle (H \cap V)$.

Applying the tidying procedure successively as in §3, given any finite subset $\mathfrak{f}$ of $\langle x_1, \ldots, x_n \rangle$ we may find a compact open subgroup, $U$, which is tidy for $\mathfrak{f}$ and such that $xUx^{-1} = U$ for every $x$ in $H \cap V$. As in Theorem 5.5, $\mathfrak{f}$ may be chosen sufficiently large that $U$ is tidy for $\langle x_1, \ldots, x_n \rangle$. Then $U$ is tidy for $\langle x_1, \ldots, x_n \rangle (H \cap V) = H$ by Proposition 5.8. □

## 6. Factoring the Scale Function and Infinitely Generated Groups of Automorphisms

Corollary 5.6 describes a decomposition of a tidy subgroup $U$ into factor subgroups $U_\mathfrak{a}$. It is shown in the present section that this decomposition is essentially unique and a number of invariants associated with it are introduced. The next definition isolates properties which characterise the factor subgroups.

**Definition 6.1.** Let $\mathfrak{H}$ be a subgroup of $\mathrm{Aut}(G)$ and $U$ be tidy for $\mathfrak{H}$. A subgroup, $V$, of $U$ satisfying:

(i). for every $\alpha$ in $\mathfrak{H}$ either $\alpha(V) \geq V$ or $\alpha(V) \leq V$;
(ii). $V = \bigcap \{\alpha(U) : \alpha \in \mathfrak{H} \text{ and } \alpha(V) \geq V\}$

will be called a *U-eigenfactor* for $\mathfrak{H}$.

The *U*-eigenfactor $\bigcap\{\alpha(U) : \alpha \in \mathfrak{H}\}$ will be denoted by $U_0$. It is a subgroup of every other *U*-eigenfactor and is invariant under all elements of $\mathfrak{H}$.

The *scale relative to the U-eigenfactor* $V$ is the function $s_V : \mathfrak{H} \to \mathbf{Z}^+$ defined by $s_V(\alpha) := [\alpha(V) : \alpha(V) \cap V]$. The set of relative scale functions on $\mathfrak{H}$ will be denoted by $S_{\mathfrak{H}}$.

The *U*-eigenfactors are precisely those groups appearing in the factoring of $U$ given in Corollary 5.6.

**Lemma 6.2.** *Let $\mathfrak{H}$ be a finitely generated subgroup of $\mathrm{Aut}(G)$ and let $U$ be tidy for $\mathfrak{H}$. Then:*

(i). *for each U-eigenfactor, $V$, there is a finite subset $\mathfrak{a}$ of $\mathfrak{H}$ such that $V = U_{\mathfrak{a}}$;*
(ii). *there are only finitely many distinct U-eigenfactors for $\mathfrak{H}$.*

**Proof.** Corollary 5.6 shows that

$$U = \prod_{\mathfrak{a} \in \Phi} U_{\mathfrak{a}},$$

where each $U_{\mathfrak{a}}$ is a *U*-eigenfactor. Now, if $\alpha(V) < V$, then $\alpha^{-1}(V) \geq V$ and so there is at least one $\mathfrak{a}$ in $\Phi$ such that $\beta(V) \geq V$ for every $\beta$ in $\mathfrak{a}$ and for this $\mathfrak{a}$ we have $V \leq U_{\mathfrak{a}}$. If $V \neq U_{\mathfrak{a}}$, then, since both groups are *U*-eigenfactors, there is an $\alpha$ in $\mathfrak{H}$ such that $\alpha(V) \geq V$ but $\alpha(U_{\mathfrak{a}}) < U_{\mathfrak{a}}$. Then $V \leq \bigcap_{k \geq 0} \alpha^k(U_{\mathfrak{a}}) = U_0$. Hence either $V = U_{\mathfrak{a}}$ for some $\mathfrak{a}$ in $\Phi$ or $V = U_0$. □

If $U_1$ and $U_2$ are tidy for $\mathfrak{H}$, with $U_1 \leq U_2$, and if $V$ is a $U_2$-eigenfactor for $\mathfrak{H}$, then $V \cap U$ is a $U_1$-eigenfactor, see Lemma 4.11. Theorem 4.12 shows that $s_{V \cap U_1} = s_V$ and so the set $S_{\mathfrak{H}}$ is independent of which tidy subgroup $U$ used to define it.

**Lemma 6.3.** *Let $U$ be a compact open subgroup of $G$ tidy for the subgroup $\mathfrak{H}$ of $\mathrm{Aut}(G)$. Let $V$ be a U-eigenfactor for $\mathfrak{H}$. Suppose that $V \neq U_0$ and put $t_V := \min\{[\alpha(V) : V] : \alpha(V) > V\}$. Choose $\alpha_V$ such that $t_V = [\alpha_V(V) : V]$. Then:*

(i). $s_V(\alpha)$ *is a power of $t_V$ for every $\alpha$ in $\mathfrak{H}$;*
(ii). *there is a surjective homomorphism $\rho_V : \mathfrak{H} \to \mathbf{Z}$, such that for every $\alpha$ in $\mathfrak{H}$*

$$\alpha(V) = \alpha_V^{\rho_V(\alpha)}(V);$$

*and*
(iii). $V = \bigcap\{\alpha(U) : \alpha \in \mathfrak{H} \text{ and } \rho_V(\alpha) \geq 0\}$.

**Proof.** (i) Let $\alpha$ and $\beta$ in $\mathfrak{H}$ be such that $\alpha(V) > V$ and $\beta(V) > V$. Then either $\alpha(V) \geq \beta(V)$, in which case $s_V(\beta)$ divides $s_V(\alpha)$, or $\alpha(V) \leq \beta(V)$, in which case $s_V(\alpha)$ divides $s_V(\beta)$. It follows from the Euclidean algorithm that, for every $\alpha$ in $\mathfrak{H}$, $s_V(\alpha)$ is a power of the greatest common divisor of $\{[\beta(V) : V] : \beta(V) > V\}$ and this, clearly, is equal to $t_V$.

(ii) The value of the homomorphism $\rho_V$ at $\alpha$ in $\mathfrak{H}$ is defined by

$$\rho_V(\alpha) = \begin{cases} \log_{t_V}[\alpha(V):V], & \text{if } \alpha(V) \geq V \\ -\log_{t_V}[V;\alpha(V)], & \text{if } \alpha(V) \leq V. \end{cases} \tag{50}$$

(iii) This follows immediately from 6.1(ii). □

It is clear, in the same way that the relative scale functions on $\mathfrak{H}$ are independent of the tidy subgroup used to define them, that the function $\rho_V$ is independent of the tidy subgroup $U$ in which $V$ is a $U$-eigenfactor. The subgroup $\mathfrak{N}_U$ defined in Theorem 4.15 is also independent of $U$ because the scale function does not depend on $U$. It is clear that this subgroup may be described in terms of the functions $\rho_V$, $V$ a $U$-eigenfactor.

**Proposition 6.4.** *Let $U$ be a compact open subgroup tidy for the finitely generated subgroup $\mathfrak{H}$ of $\mathrm{Aut}(G)$. Then*

$$\mathfrak{N} := \{\alpha \in \mathfrak{H} : \alpha(U) = U\} = \bigcap \{\ker \rho_V : V \text{ is a } U\text{-eigenfactor}\}.$$

Property **S3** of the scale function relates it to the modular function on $G$. The relative scale functions are related in the same way to modular functions on certain $\mathfrak{H}$-invariant subgroups of $G$.

**Lemma 6.5.** *Let $\mathfrak{H}$ be a group of automorphisms of $G$, let $U$ be tidy for $\mathfrak{H}$ and let $V$ be a $U$-eigenfactor. Then the group*

$$V_{++} := \bigcup_{\alpha \in \mathfrak{H}} \alpha(V)$$

*is a closed subgroup of $G$ and is invariant under $\mathfrak{H}$.*

**Proof.** Let $\alpha_V$ be the element of $\mathfrak{H}$ chosen in Lemma 6.3 and let $\mathfrak{a}$ be such that $V = U_{\mathfrak{a}}$ as in Lemma 6.2. Then $V_{++} = U_{\mathfrak{a},\alpha_V++}$ and is therefore closed, by Lemma 4.10(iv). It is clear that $V_{++}$ is invariant under $\mathfrak{H}$. □

In the case when $G$ is $\mathbf{Q}_p^n$, discussed in the introduction, the subgroups $V_{++}$ are just the common eigenspaces for the commuting automorphisms. In the general case the subgroups $V_{++}$ will be regarded as analogues of eigenspaces. It would be natural to call them eigenfactors for $\mathfrak{H}$ but they are not independent of the subgroup $U$ used to define them. The analogues of eigenvalues, which we now define, are independent of $U$ however.

**Definition 6.6.** Let $\mathfrak{H}$ be a group of automorphisms of $G$, let $U$ be tidy for $\mathfrak{H}$ and let $V$ be a $U$-eigenfactor. The function $\Delta_V : \mathfrak{H} \to \mathbf{Q}^+$ defined by

$$\Delta_V(\alpha) = \frac{m(\alpha(V))}{m(V)},$$

where $m$ is the Haar measure on $V_{++}$, is the *modulus* of $\mathfrak{H}$ relative to $V$ or the *relative modular function on* $\mathfrak{H}$. The set of relative modular functions on $\mathfrak{H}$ is denoted by $D_{\mathfrak{H}}$.

The relative modular functions on $\mathfrak{H}$ are homomorphisms and

$$\Delta_V(\alpha) = t_V^{\rho_V(\alpha)} \text{ for every } \alpha \text{ in } \mathfrak{H}. \tag{51}$$

Hence $D_{\mathfrak{H}}$ too is independent of the tidy subgroup $U$ used to define it. There is also a relative version of proper **S3** of the scale function:

$$\Delta_V(\alpha) = s_V(\alpha)/s_v(\alpha^{-1}).$$

**Definition 6.7.** Let $\mathfrak{H}$ be a finitely generated group of automorphisms of $G$ having a common tidy subgroup.

   (i). The *factor number* of $\mathfrak{H}$ is the number of distinct relative scale functions for $\mathfrak{H}$. The factor number is denoted f.n.($\mathfrak{H}$).
   (ii). The group $\mathfrak{H}/\mathfrak{N}$, where $\mathfrak{N}$ is defined in Proposition 6.4, is a finitely generated abelian group with no torsion elements by Theorem 4.15. Hence $\mathfrak{H}/\mathfrak{N} \cong \mathbf{Z}^n$. The number $n$ is the *rank* of $\mathfrak{H}$. The rank is denoted r($\mathfrak{H}$).

The factor number of $\mathfrak{H}$ is equal to the number of distinct $U$-eigenfactors for $\mathfrak{H}$ in any subgroup $U$ tidy for $\mathfrak{H}$ and also to the number of distinct non-trivial relative modular functions. Lemma 6.2 shows that f.n.($\mathfrak{H}$) is finite when $\mathfrak{H}$ is finitely generated.

Note that: f.n.($\mathfrak{H}$) = 0 means that $\alpha(U) = U$ for every $\alpha$ in $\mathfrak{H}$ and $U$ tidy for $\mathfrak{H}$; and f.n.($\mathfrak{H}$) = 1 means that either $\alpha(U) \leq U$ or $\alpha(U) \geq U$ for every $\alpha$ in $\mathfrak{H}$. Thus the common tidy subgroup $U$ is itself a $U$-eigenfactor in these cases.

**Theorem 6.8.** *Let $\mathfrak{H}$ be a finitely generated subgroup of* $\mathrm{Aut}(G)$ *and suppose that $U$ is a compact open subgroup of $G$ tidy for $\mathfrak{H}$. Then*

$$(52) \qquad U = \prod_{j=0}^{p} U_j = U_0 U_1 \ldots U_p,$$

*where $p = $ f.n.($\mathfrak{H}$) and where $U_j$, $j \in \{0, \ldots, p\}$, are distinct U-eigenfactors for $\mathfrak{H}$.*

**Proof.** Corollary 5.6 expresses $U$ as a product of $U$-eigenfactors and is seen in the proof of Lemma 6.2 that every $U$-eigenfactor appears in any such factoring. The factors need not be distinct however. To obtain distinct subgroups a sequence of automorphisms, $\beta_1, \ldots, \beta_s$ say, must be chosen more carefully.

For this, choose a fixed isomorphism $\omega : \mathbf{Z}^n \to \mathfrak{H}/\mathfrak{N}$, where $n$ is the rank of $\mathfrak{H}$ and choose elements $v_j$ in $\mathfrak{H}$ such that $\omega(e_j) = v_j\mathfrak{N}$, $j \in \{1, \ldots, n\}$, where $e_j$ denotes the $j$th unit basis vector in $\mathbf{Z}^n$. By Proposition 6.4, each homomorphism $\rho_V : \mathfrak{H} \to \mathbf{Z}$, where $V$ is a $U$-eigenfactor, induces a homomorphism $\tilde{\rho}_V : \mathfrak{H}/\mathfrak{N} \to \mathbf{Z}$. Corresponding to each of the composite maps $\tilde{\rho}_V \circ \omega : \mathbf{Z}^n \to \mathbf{Z}$, there is an $n$-tuple, $\boldsymbol{m}_V$, of integers such that

$$(53) \qquad \tilde{\rho}_V \circ \omega(\boldsymbol{x}) = \boldsymbol{m}_V.\boldsymbol{x} = \sum_{k=1}^n m_k x_k, \quad (\boldsymbol{x} \in \mathbf{Z}^n).$$

These $n$-tuples of integers will be used to choose the automorphisms $\beta_1, \ldots, \beta_s$ when there are at least two $U$-eigenfactors. Let

$$M_{\mathfrak{H}} = \{\boldsymbol{m}_V : V \text{ is a } U\text{-eigenfactor for } \mathfrak{H}\}.$$

Note that, since $\omega$ is surjective, $U$-eigenfactors $V$ and $W$ are equal if and only if $\boldsymbol{m}_V = \boldsymbol{m}_W$. Hence $|M_{\mathfrak{H}}| = $ f.n.($\mathfrak{H}$).

The cases when f.n.($\mathfrak{H}$) = 0 or 1 have been discussed above. It was seen that $U$ is itself a $U$-eigenfactor in these cases and the theorem is satisfied.

Now suppose that $|M_{\mathfrak{H}}| > 1$ and choose two points $\boldsymbol{m}_1$ and $\boldsymbol{m}_2$ in $M_{\mathfrak{H}}$. Since $M_{\mathfrak{H}}$ is finite, there is an $\boldsymbol{x}_1 = (x_1, x_2, \ldots, x_r)$ in $\mathbf{Z}^n$ such that the hyperplane
$$\Pi_{\boldsymbol{x}_1} = \{\boldsymbol{y} \in \mathbf{R}^n \mid \boldsymbol{y}.\boldsymbol{x}_1 = 0\}$$
passes between $\boldsymbol{m}_1$ and $\boldsymbol{m}_2$ and does not intersect $M_{\mathfrak{H}}$. Put
$$\beta_1 = v_1^{x_1} v_2^{x_2} \ldots v_n^{x_n}.$$
Then $\beta_1$ partitions $M_{\mathfrak{H}}$ into the two subsets
$$M_+ = \{\boldsymbol{m}_V \mid \rho_V(\beta_1) = \boldsymbol{m}_V.\boldsymbol{x}_1 > 0\} \text{ and } M_- = \{\boldsymbol{m}_V \mid \rho_V(\beta_1) = \boldsymbol{m}_V.\boldsymbol{x}_1 < 0\},$$
with, say, $\boldsymbol{m}_1$ in $M_+$ and $\boldsymbol{m}_2$ in $M_-$. Correspondingly, $S_{\mathfrak{H}}$ is partitioned into the subsets
$$S_+ = \{V \mid \beta_1(V) > V\} \text{ and } S_- = \{V \mid \beta_1(V) < V\}.$$
Now let $U = \prod_{\mathfrak{a} \in \Phi} U_{\mathfrak{a}}$ be a product of $U$-eigenfactors as in Corollary 5.6. Then by Lemma 4.4
$$U_{\beta_1+} = \prod_{\mathfrak{a} \in \Phi} (U_{\mathfrak{a}} \cap U_{\beta_1+}).$$
If $U_{\mathfrak{a}}$ is in $S_+$, then $U_{\mathfrak{a}} \cap U_{\beta_1+} = U_{\mathfrak{a}}$. On the other hand, if $U_{\mathfrak{a}}$ is in $S_-$, then $U_{\mathfrak{a}} \cap U_{\beta_1+} = U_0$ and may be absorbed into one of the other factors. Hence $U_{\beta_1+}$ is the product of subgroups in $S_+$. Similarly, $U_{\beta_1-}$ is the product of subgroups in $S_-$. Thus we have
$$U = U_{\beta_1+} U_{\beta_1-}$$
with $U_{\beta_1+}$ the product of subgroups in $S_+$ and $U_{\beta_1-}$ the product of subgroups in $S_-$. If f.n.$(\mathfrak{H}) = 2$, this shows that (52) holds.

If f.n.$(\mathfrak{H}) > 2$, then further automorphisms $\beta_2, \ldots, \beta_m$ may be chosen so that the corresponding hyperplanes separate other points in $M_{\mathfrak{H}}$. Eventually $M_{\mathfrak{H}}$ will be partitioned into singleton sets and $U$ written as a product of distinct subgroups. $\square$

It follows immediately from the above proof that distinct $U$-eigenfactors can be separated by an automorphism in $\mathfrak{H}$ which does not leave any of the $U$-eigenfactors, other than $U_0$, invariant.

**Corollary 6.9.** *Let $U$ be tidy for the group $\mathfrak{H}$ of automorphisms of $G$ and let $V_1$ and $V_2$ be $U$-eigenfactors, with $V_1 \neq V_2$ and neither equal to $U_0$. Then there is $\beta$ in $\mathfrak{H}$ such that $\beta(V_1) > V_1$, $\beta(V_2) < V_2$ and for every $U$-eigenfactor, $V$, not equal to $U_0$ we have either $\beta(V) > V$ or $\beta(V) < V$. For this $\beta$ we have $\bigcap_{k \in \mathbf{Z}} \beta^k(U) = U_0$.*

It is clear that $\mathrm{r}(\mathfrak{H}) \leq \mathrm{f.n.}(\mathfrak{H})$. It is also clear that, if $\mathrm{r}(\mathfrak{H}) = 0$, then $\mathrm{f.n.}(\mathfrak{H}) = 0$ and, if $\mathrm{r}(\mathfrak{H}) = 1$, then $\mathrm{f.n.}(\mathfrak{H}) = 1$ or 2. These are the only relationships between the rank and the factor number of $\mathfrak{H}$ however. The following example shows that, if $\mathrm{r}(\mathfrak{H}) = 2$, then $\mathrm{f.n.}(\mathfrak{H})$ can take any value except 0 or 1 and the same method produces similar examples with $\mathrm{r}(\mathfrak{H}) > 2$.

**Example 6.10.** Let $\Psi$ be a finite set in $\mathbf{Z}^2$ which contains $\boldsymbol{0}$ and is such that $\gcd\{m_1, m_2\} = 1$ for each $\boldsymbol{m} = (m_1, m_2)$ in $\Psi$. Define the group
$$G = \prod_{\boldsymbol{m} \in \Psi} \mathbf{Q}_3$$
and write the elements of $G$ as functions $g : \Psi \to \mathbf{Q}_3$. Define autmorphisms of $G$
$$(\alpha_1 g)(\boldsymbol{m}) = 3^{-m_1} g(\boldsymbol{m}) \text{ and } (\alpha_2 g)(\boldsymbol{m}) = 3^{-m_2} g(\boldsymbol{m})$$

and let $\mathfrak{H}$ be the subgroup of $\mathrm{Aut}(G)$ generated by $\alpha_1$ and $\alpha_2$. Then $U := \prod_{\boldsymbol{m} \in \Psi} \mathbf{Z}_3$ is tidy for $\mathfrak{H}$ and the $U$-eigenfactors for $\mathfrak{H}$ are

$$U_{\boldsymbol{m}} = \{g \in U \mid g(\boldsymbol{m}') = 0 \text{ unless } \boldsymbol{m}' = \boldsymbol{m} \text{ or } \boldsymbol{0}\}.$$

Now consider an $\boldsymbol{m}$ in $\Psi$. Since $\gcd\{m_1, m_2\} = 1$, there is an $\alpha_0$ in $\mathfrak{H}$ such that $\alpha_0 g = 3^{-1} g$ for $g$ in $U_{\boldsymbol{m}}$. Hence the number $t_{\boldsymbol{m}}$ in Lemma 6.3 is 3 and the homomorphism $\rho_{\boldsymbol{m}} : \mathfrak{H} \to \mathbf{Z}$ determined by $U_{\boldsymbol{m}}$ is

$$\rho_{\boldsymbol{m}}(\alpha) = m_1 x_1 + m_2 x_2, \quad (\alpha = \alpha_1^{x_1} \alpha_2^{x_2} \in \mathfrak{H}).$$

Therefore, choosing the isomorphism $\mathfrak{H} \to \mathbf{Z}^2$ determined by $\alpha_1 \mapsto (1,0)$ and $\alpha_2 \mapsto (0,1)$, we find that $M_{\mathfrak{H}} = \Psi \setminus \{\boldsymbol{0}\}$.

In the previous example the tidy subgroup $U$ is abelian and the $U$-eigenfactors in the product (52) may be written in any order. That is not so in general.

**Example 6.11.** Let $G = SL(3, \mathbf{Q}_p)$ and $\mathfrak{H} = \{\alpha_d : d \in G \text{ is diagonal}\}$ be the group of inner automorphisms of $G$ determined by the diagonal matrices. Then the compact open subgroup

$$U := \{a = [a_{ij}] \in G : a_{ij} \in \mathbf{Z}_p \text{ if } i \geq j, \ a_{ij} \in p\mathbf{Z}_p \text{ if } i < j\}$$

is tidy for $\mathfrak{H}$, [3]. It is clear that $\mathfrak{N}$ is the subgroup of $\mathfrak{H}$ corresponding to the diagonal matrices with $p$-adic integer entries and that $\mathfrak{H}/\mathfrak{N} \cong \mathbf{Z}^2$. Hence $\mathrm{r}(\mathfrak{H}) = 2$.

Let $d^{(1)} = \mathrm{diag}[p^{-1}, 1, p]$ and let $\alpha_1$ be the automorphism of conjugation by $d^{(1)}$. Then

$$U_{\alpha_1+} = \left\{ \begin{bmatrix} a_{11} & a_{12} & a_{13} \\ 0 & a_{22} & a_{23} \\ 0 & 0 & a_{33} \end{bmatrix} : a_{ii} \in \mathbf{Z}_p, \ a_{ij} \in p\mathbf{Z}_p, \ i < j \right\}$$

$$\text{and} \quad U_{\alpha_1-} = \left\{ \begin{bmatrix} a_{11} & 0 & 0 \\ a_{21} & a_{22} & 0 \\ a_{31} & a_{32} & a_{33} \end{bmatrix} : a_{ij} \in \mathbf{Z}_p \right\}$$

and property $\mathbf{T1}(\alpha_1)$ expresses the $UL$-factorisation of matrices in $U$.

Let $d^{(2)} = \mathrm{diag}[p^{-1}, p, 1]$ and $d^{(3)} = \mathrm{diag}[1, p^{-1}, p]$ and $\alpha_2$ and $\alpha_3$ be the corresponding inner automorphisms. Then, taking $\mathfrak{a} = \{\alpha_1, \alpha_2^{-1}\}$, the subgroup

$$U_{\mathfrak{a}} = \left\{ \begin{bmatrix} a_{11} & 0 & 0 \\ 0 & a_{22} & a_{23} \\ 0 & 0 & a_{33} \end{bmatrix} : a_{ii} \in \mathbf{Z}_p, \ a_{23} \in p\mathbf{Z}_p \right\}$$

is a $U$-eigenfactor. It is also a root subgroup of $U$. The relative modular function corresponding to this $U$-eigenfactor is $\Delta_{\mathfrak{a}}(\alpha_d) = |d_2 d_3^{-1}|_p$ and the relative scale function is $s_{\mathfrak{a}}(\alpha_d) = \max\{|d_2 d_3^{-1}|_p, 1\}$, where $d = \mathrm{diag}[d_1, d_2, d_3]$.

The groups $U_{\mathfrak{b}}$ and $U_{\mathfrak{c}}$, where $\mathfrak{b} = \{\alpha_1, \alpha_3^{-1}\}$ and $\mathfrak{c} = \{\alpha_1, \alpha_2, \alpha_3\}$, are also $U$-eigenfactors and root subgroups of $U$ and it is easily checked that $U_{\alpha_1+} = U_{\mathfrak{a}} U_{\mathfrak{c}} U_{\mathfrak{b}}$. Similarly, $U_{\alpha_1-}$ is a product of three $U$-eigenfactors. Hence $U$ is the product of six $U$-eigenfactors (other than $U_0$) and $\mathrm{f.n.}(\mathfrak{H}) = 6$.

The factoring of subgroups tidy for $\mathfrak{H}$ allows us to write the scale and modular functions (restricted to $\mathfrak{H}$) as products of relative scale and modular functions.

**Theorem 6.12.** *Let $\mathfrak{H}$ be a finitely generated subgroup of* $\mathrm{Aut}(G)$ *and suppose that $\mathfrak{H}$ has a tidy subgroup. Then for each $\alpha$ in $\mathfrak{H}$*

$$s(\alpha) = \prod \{s_V(\alpha) : s_V \in S_{\mathfrak{H}}\} \quad and \quad \Delta(\alpha) = \prod \{\Delta_V(\alpha) : \Delta_V \in D_{\mathfrak{H}}\}.$$

**Proof.** Let $U$ be tidy for $\mathfrak{H}$. Then $U = U_{\alpha+} U_{\alpha-}$, where $U_{\alpha+} = \prod_{j=0}^{p} U_j$ is the product of those $U$-eigenfactors such that $\alpha(U_j) \geq U_j$. If $p = 1$, then $s(\alpha) = [\alpha(U_1) : U_1] = s_{U_1}(\alpha)$ and the first equation is proved.

Suppose that $p > 1$. Then, by Corollary 6.9, there is $\beta$ in $\mathfrak{H}$ such that $\beta(U_1) > U_1$, $\beta(U_2) < U_2$ and $U_{\alpha+,\beta 0} = U_0$. Then equation (35) in Lemma 5.2 shows that

$$s(\alpha) = [\alpha(U_{\alpha+}) : U_{\alpha+}] = [\alpha(U_{\alpha+,\beta+}) : U_{\alpha+,\beta+}][\alpha(U_{\alpha+,\beta-}) : U_{\alpha+,\beta-}].$$

Now $U_{\alpha+,\beta-}$ and $U_{\alpha+,\beta+}$ are each the product of $U$-eigenfactors on which $\alpha$ is non-contracting and, by the choice of $\beta$, $U_{\alpha+,\beta-} \cap U_{\alpha+,\beta+} = U_0$. If $U_{\alpha+,\beta-}$ and $U_{\alpha+,\beta+}$ are both $U$-eigenfactors, then the formula for $s(\alpha)$ is proved. Otherwise, further automorphisms, $\beta'$, may be chosen with the aid of Corollary 6.9 such that $[\alpha(U_{\alpha+,\beta+}) : U_{\alpha+,\beta+}]$ and $[\alpha(U_{\alpha+,\beta-}) : U_{\alpha+,\beta-}]$ are factored. Continuing in this way, $s(\alpha)$ is eventually written as the product of relative scale functions.

That $\Delta(\alpha)$ is the product of relative modular functions follows from the formula for $s(\alpha)$ because $\Delta(\alpha) = s(\alpha)/s(\alpha^{-1})$. □

So far in this section we have dealt mainly with finitely generated groups of automorphisms having a tidy subgroup. The definitions and results will now be extended to infinitely generated groups of automorphisms having local tidy subgroups.

**Lemma 6.13.** *Let $\mathfrak{H}$ be a subgroup of* $\mathrm{Aut}(G)$ *having local tidy subgroups. Let $\mathfrak{F}$ be a finitely generated subgroup of $\mathfrak{H}$.*

(i). *Then there is an integer $r$ and there are (not necessarily distinct) homomorphisms $\phi_j : \mathfrak{F} \to (\mathbf{Q}^+, \times)$, $j \in \{1, \ldots, r\}$ such that, for every finitely generated group $\mathfrak{K}$ with $\mathfrak{F} \leq \mathfrak{K} \leq \mathfrak{H}$, there are mutually disjoint subsets $S_\Delta$ of $\{1, \ldots, r\}$ with*

$$\Delta(\alpha) = \prod_{j \in S_\Delta} \phi_j(\alpha), \quad (\alpha \in \mathfrak{F}, \Delta \in D_{\mathfrak{K}}),$$

*and $\bigcup_{\Delta \in D_{\mathfrak{K}}} S_\Delta = \{1, \ldots, r\}$.*

(ii). *There is a finitely generated group $\mathfrak{L}$ with $\mathfrak{F} \leq \mathfrak{L} \leq \mathfrak{H}$ such that for each $j$ in $\{1, \ldots, r\}$ there is a relative modular function $\Delta^{(j)}$ in $D_{\mathfrak{L}}$ with $\phi_j = \Delta^{(j)}|_{\mathfrak{F}}$ and with $\Delta^{(i)} \neq \Delta^{(j)}$ if $i \neq j$.*

(Note in (i) that, if $\Delta|_{\mathfrak{F}}$ is the trivial homomorphism for some $\Delta$ in $D_{\mathfrak{K}}$, then $S_\Delta = \emptyset$ and so $\{S_\Delta : \Delta \in D_{\mathfrak{K}}\}$ may not be a partition of $\{1, \ldots, r\}$.)

**Proof.** Let $\mathfrak{F}'$ be a finitely generated subgroup of $\mathfrak{H}$ with $\mathfrak{F}' \geq \mathfrak{F}$. Let $U$ be a compact open subgroup of $G$ which is tidy for $\mathfrak{F}'$. Then $U$ is also tidy for $\mathfrak{F}$ and so $U$ may factored as $U = \prod_{j=0}^{p} U_j$, where $U_j$, $j \in \{1, \ldots, p\}$ are $U$-eigenfactors for $\mathfrak{F}$ and this may be refined to $U = \prod_{k=0}^{q} U'_k$ where $U'_k$, $k \in \{1, \ldots, q\}$ are $U$-eigenfactors for $\mathfrak{F}'$. Let $\Delta_j$, $j \in \{1, \ldots, p\}$ (respectively $\Delta'_k$, $k \in \{1 \ldots, q\}$) be the corresponding relative modular functions on $\mathfrak{F}$ (respectively $\mathfrak{F}'$). That the second factoring of $U$ is a refinement of the first means that each $U$-eigenfactor, $U_j$, for $\mathfrak{F}$ is a product $U_j = \prod_{k=k_j}^{k_{j+1}-1} U'_k$. It now follows as in Theorem 6.12 that for each $j \in \{1,\ldots,p\}$ and each $\alpha$ in $\mathfrak{F}$ which is expanding on $U_j$

$$(54) \qquad \Delta_j(\alpha) = [\alpha(U_j) : U_j] = \prod_{k=k_j}^{k_{j+1}-1} [\alpha(U'_k) : U'_k] = \prod_{k=k_j}^{k_{j+1}-1} \Delta'_k(\alpha).$$

It is clear that the factoring of $\Delta_j$ in (54) is further refined when another finitely generated subgroup, $\mathfrak{F}''$, of $\mathfrak{H}$ with $\mathfrak{F}'' \geq \mathfrak{F}'$ is taken. Now the values of all relative modular functions are rational and consequently have only finitely many factors. Hence there is a finitely generated subgroup, $\mathfrak{L}$, of $\mathfrak{H}$ such that (54) holds with $\mathfrak{F}' = \mathfrak{L}$ and for every finitely generated $\mathfrak{F}'' \geq \mathfrak{L}$ the corresponding factorisation $\Delta_j(\alpha) = \prod_{l=l_j}^{l_{j+1}-1} \Delta''_l(\alpha)$, where $\Delta''_l$ in $D_{\mathfrak{F}''}$, satisfies

$$(55) \qquad \text{either } \Delta''_l|_{\mathfrak{F}} = \Delta'_k|_{\mathfrak{F}} \text{ for some } k \text{ or } \Delta''_l|_{\mathfrak{F}} \text{ is trivial.}$$

Each relative modular function, $\Delta'$, on $\mathfrak{L}$ restricts to a homomorphism from $\mathfrak{F}$ to $(\mathbf{Q}^+, \times)$. Let $r$ be the number of distinct relative modular functions on $\mathfrak{L}$ which have non-trivial, but not necessarily distinct, restrictions and define $\phi_j = \Delta'|_{\mathfrak{F}}$, $j \in \{1,\ldots,r\}$, to be these restrictions. Then (ii) is satisfied.

To see that these functions satisfy (i), let $\mathfrak{K} \geq \mathfrak{F}$ be a finitely generated subgroup of $\mathfrak{H}$ and put $\mathfrak{F}'' = \langle \mathfrak{K}, \mathfrak{L} \rangle$. Then (54) shows that each relative modular function, $\Delta$, on $\mathfrak{K}$ is the product of relative modular functions on $\mathfrak{F}''$. Hence $\Delta|_{\mathfrak{F}}$ is the product of restrictions to $\mathfrak{F}$ of functions in $D_{\mathfrak{F}''}$. It now follows from (55) that $\Delta|_{\mathfrak{F}}$ is in fact the product of functions $\phi_j$. If $\Delta_1$ and $\Delta_2$ are distinct elements of $D_{\mathfrak{K}}$, then they are the product of disjoint sets of relative modular functions on $\mathfrak{F}''$. Now distinct elements of $\mathfrak{F}''$ restrict to the trivial function or to functions, $\phi_i$ and $\phi_j$ with $i \neq j$. Hence the sets $S_{\Delta_1}$ and $S_{\Delta_2}$ are disjoint. $\square$

**Theorem 6.14.** *Let $\mathfrak{H}$ be a subgroup of* $\mathrm{Aut}(G)$ *having local tidy subgroups. Then there is a set $\Phi \subseteq \mathrm{Hom}(\mathfrak{H}, (\mathbf{Q}^+, \times))$, such that:*

(i). *for each $\phi \in \Phi$ there are: a positive integer $t_\phi$; and a homomorphism $\rho_\phi : \mathfrak{H} \to \mathbf{Z}$ such that $\phi(\alpha) = t_\phi^{\rho_\phi(\alpha)}$ for every $\alpha \in \mathfrak{H}$;*

(ii). $\bigcap_{\phi \in \Phi} \ker \phi = \{\alpha \in \mathfrak{H} : s(\alpha) = 1 = s(\alpha^{-1})\} := \mathfrak{N}$;

(iii). $\{\phi \in \Phi : \phi(\alpha) \neq 1\}$ *is finite for each $\alpha$ in $\mathfrak{H}$; and*

(iv). *if $\mathfrak{F} \leq \mathfrak{H}$ is finitely generated and $\Delta$ is in $D_{\mathfrak{F}}$, then $\Delta(\alpha) = \prod_{\phi \in \Phi} \phi(\alpha)$ for every $\alpha$ in $\mathfrak{F}$.*

**Proof.** For each finitely generated subgroup, $\mathfrak{F}$, of $\mathfrak{H}$ the homomorphisms $\phi_j : \mathfrak{F} \to (\mathbf{Q}^+, \times)$, $j \in \{1,\ldots,r\}$ found in Lemma 6.13 might not be distinct. Form the multiset, $\Phi_{\mathfrak{F}}$, of these homomorphisms, so that $\phi_i$ and $\phi_j$ are distinct elements of $\Phi_{\mathfrak{F}}$ if $i \neq j$ even though they might be the same homomorphism.

Suppose that $\mathfrak{F} \leq \mathfrak{F}'$ are finitely generated subgroups of $\mathfrak{H}$. Let $\mathfrak{L} \geq \mathfrak{F}'$ be the finitely generated subgroup of $\mathfrak{H}$ found in Lemma 6.13(ii). Then for every $\phi_j$ in $\Phi_{\mathfrak{F}'}$ there is a unique $\Delta^{(j)}$ in $D_{\mathfrak{L}}$ such that $\Delta^{(j)}|_{\mathfrak{F}'} = \phi_j$. On the other hand, if $\Delta$ is in $D_{\mathfrak{L}}$, then either $\Delta|_{\mathfrak{F}}$ is trivial or $\Delta|_{\mathfrak{F}} = \phi_j$ for some $j$ and, if $\Delta_1|_{\mathfrak{F}} = \phi_j = \Delta_2|_{\mathfrak{F}}$, then $\Delta_1 = \Delta_2$. Since $\mathfrak{F} \leq \mathfrak{F}'$, the same statements hold with $\mathfrak{F}$ in place of $\mathfrak{F}'$. It follows that for each $\phi$ in $\Phi_{\mathfrak{F}}$ there is a unique $\phi'$ in $\Phi_{\mathfrak{F}'}$ which extends it. Since this holds for every $\mathfrak{F}' \geq \mathfrak{F}$, each $\phi$ in $\Phi_{\mathfrak{F}}$ has a unique extension, $\hat{\phi}$, to $\mathfrak{H}$ such that

$\hat{\phi}|_{\mathfrak{F}'}$ is in $\Phi_{\mathfrak{F}'}$ for every $\mathfrak{F}' \geq \mathfrak{F}$. Put

$$\Phi = \{\hat{\phi} : \phi \in \Phi_{\mathfrak{F}}, \mathfrak{F} \leq \mathfrak{H} \text{ and finitely generated}\}.$$

Note that distinct elements of $\Phi_{\mathfrak{F}}$ have distinct extensions to $\mathfrak{L}$, by Lemma 6.13(ii), and so the map $\phi \mapsto \hat{\phi}$, from the multiset $\Phi_{\mathfrak{F}}$ to the set $\Phi$ is injective.

The existence of the positive integers $t_\phi$ and homomorphisms $\rho_\phi : \mathfrak{H} \to \mathbf{Z}$ now follow from Lemma 6.3(ii). Let $\alpha$ be in $\mathfrak{H}$, put $\mathfrak{F} = \langle \alpha \rangle$ and let $\mathfrak{L}$ be the subgroup of $\mathfrak{H}$ found in Lemma 6.13(ii) for this $\mathfrak{F}$. Then it follows from Theorem 6.12 that $s(\alpha) = \prod \{\phi(\alpha) : \phi \in \Phi_{\mathfrak{L}}, \phi(\alpha) \geq 1\}$. Hence, if $\alpha$ belongs to $\bigcap_{\phi \in \Phi} \ker \phi$, then $s(\alpha) = 1$ and similarly for $s(\alpha^{-1})$. On the other hand, if $s(\alpha) = 1 = s(\alpha^{-1})$, then $\phi(\alpha) = 1$ for every $\phi$ in $\Phi_{\mathfrak{L}}$, whence $\hat{\phi}(\alpha) = 1$ for every $\phi$ in $\Phi_{\mathfrak{L}}$. Since $\phi(\alpha) = 1$ for every $\phi$ in $\Phi \setminus \hat{\Phi}_{\mathfrak{L}}$, it follows that $\alpha$ belongs to $\bigcap_{\phi \in \Phi} \ker \phi$. Thus (ii) is established. It is clear that $\phi(\alpha) \neq 1$ only if $\phi$ is in $\hat{\Phi}_{\mathfrak{L}}$ and so (iii) holds. Part (iv) follows from Theorem 6.12. □

It follows from parts (i) and (ii) of the Theorem that the map

(56) $$R := \bigoplus_{\phi \in \Phi} \rho_\phi : \mathfrak{H} \to \bigoplus_{\phi \in \Phi} \mathbf{Z}$$

induces an isomorphism between $\mathfrak{H}/\mathfrak{N}$ and a subgoup of the free abelian group $\bigoplus_{\phi \in \Phi} \mathbf{Z}$. Since each subgroup of a free abelian group is free, we have the following improvement on Theorem 4.15.

**Corollary 6.15.** *Let $\mathfrak{H}$ be a subgroup of $\mathrm{Aut}(G)$ having local tidy subgroups and define $\mathfrak{N} = \{\alpha \in \mathfrak{H} : s(\alpha) = 1 = s(\alpha^{-1})\}$. Then $\mathfrak{N}$ is normal in $\mathfrak{H}$ and $\mathfrak{H}/\mathfrak{N}$ is a free abelian group.*

The quotient group $A_{\mathfrak{H}} := \left(\bigoplus_{\phi \in \Phi} \mathbf{Z}\right)/R(\mathfrak{H})$, where $R$ is the map defined in (56), need be not be free abelian. In fact, any abelian group can occur as such a quotient, as we now see.

Let $\mathfrak{c}$ be a cardinal number, $p$ be a prime, and let $G$ be the restricted product of copies of $\mathbf{Q}_p$,

$$G = \{g \in \mathbf{Q}_p^{\mathfrak{c}} : g(c) \in \mathbf{Z}_p \text{ for all but finitely many } c \in \mathfrak{c}\}.$$

For each $\boldsymbol{h}$ in $\bigoplus_{\mathfrak{c}} \mathbf{Z}$ let $U_{\boldsymbol{h}}$ be the subgroup $\{g \in G : g(c) \in p^{\boldsymbol{h}(c)}\mathbf{Z}_p \text{ for all } c \in \mathfrak{c}\}$. Then $\{U_{\boldsymbol{h}}\}$ is a base of neighbourhoods of the identity for a group topology on $G$. In this topology $U_{\boldsymbol{0}}$ is isomorphic to $\mathbf{Z}^{\mathfrak{c}}$ with the product topology and is a compact open subgroup of $G$. Each $\boldsymbol{h}$ in $\bigoplus_{\mathfrak{c}} \mathbf{Z}$ determines a continuous automorphism of $G$ by

$$(\boldsymbol{h}g)(c) = p^{-\boldsymbol{h}(c)}g(c), \quad (c \in \mathfrak{c}).$$

This automorphism will also be denoted by $\boldsymbol{h}$ and the group of these automorphisms identified with $\bigoplus_{\mathfrak{c}} \mathbf{Z}$. Then $U_{\boldsymbol{0}}$ is tidy for $\bigoplus_{\mathfrak{c}} \mathbf{Z}$. Let $\mathfrak{H}$ be a subgroup of $\bigoplus_{\mathfrak{c}} \mathbf{Z}$ and suppose that

(a) $\gcd\{\boldsymbol{h}(c) : \boldsymbol{h} \in \mathfrak{H}\} = 1, \quad (c \in \mathfrak{c})$
(b) and for $c \neq d \in \mathfrak{c}$, there is $\boldsymbol{h} \in \mathfrak{H}$ with $\boldsymbol{h}(c) > 0$ and $\boldsymbol{h}(d) < 0$.

Condition (a) implies that for each $c$ in $\mathfrak{c}$ there is an $\boldsymbol{h}$ in $\mathfrak{H}$ such that $\boldsymbol{h}(c) = 1$. Then $\boldsymbol{h}$ has finite support in $\mathfrak{c}$ and condition (b) implies that, for each $d \neq c$ in the support of $\boldsymbol{h}$, there is $\boldsymbol{h}_d$ such that $\boldsymbol{h}_d(d) < 0$ and $\boldsymbol{h}_d(c) > 0$. Put

$$\mathfrak{F} = \langle \boldsymbol{h}, \boldsymbol{h}_d : d \in \mathrm{supp}(\boldsymbol{h}),\ d \neq c \rangle$$

and $\mathfrak{h} = \{\boldsymbol{h}_d : d \in \mathrm{supp}(\boldsymbol{h}),\ d \neq c\}$. Then

$$U_\mathfrak{h} = \bigcap \{\boldsymbol{h}_d^n(U) : d \in \mathrm{supp}(\boldsymbol{h}),\ d \neq c,\ n \geq 0\}$$

is a compact subgroup of $U_\mathbf{0}$ and $[\boldsymbol{h}(U_\mathfrak{h}) : U_\mathfrak{h}] = p$. Now $U_\mathfrak{h}$ might not be a $U_\mathbf{0}$-eigenfactor for $\mathfrak{F}$ but, adding further elements from $\mathfrak{F}$ to $\mathfrak{h}$ if necessary, a subgroup, $U'$, may be found which is a $U_\mathbf{0}$-eigenfactor for $\mathfrak{F}$ and such that $[\boldsymbol{h}(U') : U'] = p$. Hence the homomorphism $\phi_c : \boldsymbol{h} \mapsto p^{-\boldsymbol{h}(c)}$ belongs to the set $\Phi$ defined in Theorem 6.14. It is clear then that $\Phi = \{\phi_c\}_{c \in \mathfrak{c}}$. Under the identification we have made between $\bigoplus_\mathfrak{c} \mathbf{Z}$ and a subgroup of $\mathrm{Aut}(G)$, it is clear that $R(\mathfrak{H}) = \mathfrak{H}$. Therefore

$$A_\mathfrak{H} \cong \left( \bigoplus_\mathfrak{c} \mathbf{Z} \right) / \mathfrak{H}.$$

Finally, we show that each abelian subgroup is isomorphic to a group $\left( \bigoplus_\mathfrak{c} \mathbf{Z} \right) / \mathfrak{H}$ where $\mathfrak{H}$ satisfies (a) and (b). Let $A$ be an abelian group. Then $A \cong \left( \bigoplus_\mathfrak{c} \mathbf{Z} \right) / \mathfrak{H}$ for some free abelian group $\bigoplus_\mathfrak{c} \mathbf{Z}$ and subgroup $\mathfrak{H}$.

Suppose that $\gcd\{\boldsymbol{h}(c) : \boldsymbol{h} \in \mathfrak{H}\} = m \neq 1$ for some $c$ in $\mathfrak{c}$. In the case when $m \neq 0$, let $\eta : \mathbf{Z}^m \to \mathbf{Z}$ be the homomorphism $\eta(x_1, \ldots, x_m) = x_1 + \cdots + x_m$, let $\phi_c : \bigoplus_\mathfrak{c} \mathbf{Z} \to \mathbf{Z}$ be the homomorphism $\phi_c(\boldsymbol{h}) = \boldsymbol{h}(c)$ and let $\bigoplus_\mathfrak{c} \mathbf{Z} \times_\mathbf{Z} \mathbf{Z}^m$ be the pull-back which completes the commutative diagram

$$\begin{array}{ccc} \bigoplus_\mathfrak{c} \mathbf{Z} \times_\mathbf{Z} \mathbf{Z}^m & \xrightarrow{\psi_c} & \mathbf{Z}^m \\ \varsigma \downarrow & & \downarrow \eta \\ \bigoplus_\mathfrak{c} \mathbf{Z} & \xrightarrow{\phi_c} & \mathbf{Z} \end{array}$$

Then $\bigoplus_\mathfrak{c} \mathbf{Z} \times_\mathbf{Z} \mathbf{Z}^m$ is a free abelian group and in fact there is an obvious isomorphism $\bigoplus_\mathfrak{c} \mathbf{Z} \times_\mathbf{Z} \mathbf{Z}^m \cong \bigoplus_{\mathfrak{c} \setminus \{c\}} \mathbf{Z} \times \mathbf{Z}^m$ which replaces the $c$-component of $\bigoplus_\mathfrak{c} \mathbf{Z}$ by $\mathbf{Z}^m$. Define $\mathfrak{H}'$ to be the subgroup $\varsigma^{-1}(\mathfrak{H})$. Then $(\bigoplus_\mathfrak{c} \mathbf{Z} \times_\mathbf{Z} \mathbf{Z}^m)/\mathfrak{H}' \cong A$. It may be seen that

$$\eta^{-1}(\phi_c(\mathfrak{H})) = \mathbf{Z}(1, \ldots, 1) \oplus \ker \eta = \langle (1, \ldots, 1), (1, -1, 0, \ldots), \ldots, (1, 0, \ldots, 0, -1) \rangle$$

and this group satisfies (a). It follows that $\mathfrak{H}' = \psi_c^{-1}(\eta^{-1}(\phi_c(\mathfrak{H})))$ satisfies (a) as a subgroup of $\bigoplus_{\mathfrak{c} \setminus \{c\}} \mathbf{Z} \times \mathbf{Z}^m$. The case when $m = 0$ can be treated by a similar argument using the map $\eta : \mathbf{Z}^2 \to \mathbf{Z}$, $\eta(x_1, x_2) = x_1 + x_2$ and noting that in this case $\phi_c(\mathfrak{H}) = \{0\}$. Since the $\bigoplus_\mathfrak{c} \mathbf{Z}$ is the direct sum, repeating this argument for each $c$ in $\mathfrak{c}$ will produce a group satisfying (a).

Suppose that there are distinct $c_1$ and $c_2$ in $\mathfrak{c}$ such that $\boldsymbol{h}(c_1)$ and $\boldsymbol{h}(c_2)$ have the same sign for every $\boldsymbol{h}$ in $\mathfrak{H}$. Then, since $\mathfrak{H}$ satisfies (a), it follows that in fact $\boldsymbol{h}(c_1) = \boldsymbol{h}(c_2)$ for every $\boldsymbol{h}$ in $\mathfrak{H}$. Put $\mathfrak{f} = \{d \in \mathfrak{c} : \boldsymbol{h}(d) = \boldsymbol{h}(c_1)\ \forall \boldsymbol{h} \in \mathfrak{H}\}$ and define an automorphism, $\nu$, of $\bigoplus_\mathfrak{c} \mathbf{Z}$ by

$$\nu \boldsymbol{h}(d) = \begin{cases} \boldsymbol{h}(d) & \text{if } d \in \mathfrak{c} \setminus \mathfrak{f}\ \text{or}\ d = c_1 \\ \boldsymbol{h}(d) - \boldsymbol{h}(c_1) & \text{if } d \in \mathfrak{f} \setminus \{c_1\}. \end{cases}$$

Then $(\bigoplus_{\mathfrak{c}} \mathbf{Z})/\nu(\mathfrak{H}) \cong A$ but $\boldsymbol{h}(d) = 0$ for every $\boldsymbol{h}$ in $\nu(\mathfrak{H})$ and every $d$ in $\mathfrak{f} \setminus \{c_1\}$, so that $\nu(\mathfrak{H})$ does not satisfy (a). However, applying the $m = 0$ case of previous argument to each element of $\mathfrak{f} \setminus \{c_1\}$ ensures that (a) is satisfied. Repeating this argument for each level set of $\mathfrak{H}$ ensures that (b) is satisfied as well as (a).

**Definition 6.16.** Let $\mathfrak{H}$ be a subgroup of $\mathrm{Aut}(G)$ having local tidy subgroups.
  (i). Let $\Phi$ be the set of homomorphisms $\mathfrak{H} \to \mathbf{Z}$ defined in Theorem 6.14. Then $|\Phi|$, the cardinality of $\Phi$, is the *factor number* of $\mathfrak{H}$ and is denoted f.n.($\mathfrak{H}$).
  (ii). Let $R$ be the homomorphism $\mathfrak{H} \to \bigoplus_\Phi \mathbf{Z}$ defined in (56). Then $R(\mathfrak{H}) \cong \mathfrak{H}/\mathfrak{N}$ is a free abelian group and is the *rank group* of $\mathfrak{H}$. The rank of $R(\mathfrak{H})$ is the *rank* of $\mathfrak{H}$ and is denoted r($\mathfrak{H}$).
  (iii). The group $(\bigoplus_\Phi \mathbf{Z})/R(\mathfrak{H})$ is the *corank group* of $\mathfrak{H}$ and denoted $A_\mathfrak{H}$.

Let $G = SL(3, \mathbf{Q}_p)$ and $\mathfrak{H}$ be as in Example 6.11. Then r($\mathfrak{H}$) = 2, f.n.($\mathfrak{H}$) = 6 and it may be seen that $A_\mathfrak{H} \cong \mathbf{Z}^4$. The group $\mathfrak{H}$ of inner automorphisms in this example is a maximal group for which there is a tidy subgroup and the rank of $\mathfrak{H}$ is the same as the usual rank of $G$. It might be thought that the rank of groups of inner automorphisms maximal with respect to the property of having local tidy subgroups could be used to extend the definition of rank to arbitrary totally disconnected locally compact groups. Any such extension of the notion of rank could not be straightforward however, as the following example shows.

**Example 6.17.** Let $E = \langle c_1, c_2, a : c_i^2 = e_E = a^2, c_1 c_2 = c_2 c_1, a c_1 = c_2 a \rangle$ and define $C = \langle c_1 \rangle = \{e_E, c_1\} \leq E$. Let
$$H = \{h \in E^\mathbf{Z} : h(n) \in C \text{ for almost all } n \in \mathbf{Z}\}.$$
For each $n \geq 0$ let $U_n$ be the subgroup of $H$,
$$U_n = \{h \in C^\mathbf{Z} : h(j) = e_E \text{ for } -n < j < n\}$$
and topologise $H$ by letting $\{U_n\}_{n \geq 0}$ be a base of neighbourhoods of the identity. Define the shift automorphism $\alpha$ of $H$ by
$$(\alpha h)(n) = h(n+1), \quad (h \in H)$$
and let $G = H \times_\alpha \mathbf{Z}$ be the semi-direct product.

The normaliser of $U_0$ is the group
$$B := \{h \in H : h(n) \in \langle c_1, c_2 \rangle\} \times_\alpha \mathbf{Z}$$
and $U_0$ is the only compact open subgroup tidy for $B$. If $g$ belongs to $G \setminus B$, then $gU_0 g^{-1}$ is not tidy for $B$ and so, by Lemma 4.3, $U_0$ is not tidy for $\langle g, A \rangle$. Hence $B$ is a maximal group among subgroups of $B$ having a tidy subgroup. Since $U_0$ is normal in $B$, f.n.($B$) = r($B$) = 0 and $A_B$ is the trivial group.

Let $k$ be the element of $H$ given by $k(n) = e_E$ if $n \neq 0$ and $k(0) = a$ and let $\beta = k\alpha$. Then conjugation by $\beta$ is an automorphism of $H$. Define
$$D := \{h \in H : h(n) \in \langle c_1, c_2 \rangle\} \times_\beta \mathbf{Z} < G.$$
The compact open subgroup $U_0$ is tidy for $D$ and $D$ is maximal among subgroups of $G$ having a tidy subgroup. Now $U_0$ is normal in $\{h \in H : h(n) \in \langle c_1, c_2 \rangle\}$ but
$$\beta U_0 \beta^{-1} = \{h \in H : h(j) = \{e_E, c_1\} \text{ if } j \neq -1 \text{ and } h(-1) \in \{e_E, c_2\}\}.$$

Hence the homomorphisms from $D$ to $\mathbf{Q}^+$ in the set $\Phi$ defined in Theorem 6.14 are determined by their values on $\beta$. There are just two of them: $\phi_1(\beta) = 2$ and $\phi_2(\beta) = \frac{1}{2}$. The map $R : D \to \mathbf{Z}^2$ is given by

$$R(h\beta^n) = (n, -n), \quad (h\beta^n \in D).$$

Therefore f.n.$(D) = 2$, r$(D) = 1$ and $A_D \cong \mathbf{Z}$.

The rank and factor number of those subgroups of $G$ which are maximal with respect to having a tidy subgroup are thus not unique.

In addition to the notion of rank for a group of automorphisms, $\mathfrak{H}$, having local tidy subgroups, there is an analogue of the Weyl group. The following theorem shows that the normaliser of $\mathfrak{H}$ in $\mathrm{Aut}(G)$ has a representation on the free abelian group $\bigoplus_{\phi \in \Phi} \mathbf{Z}$ and this representation induces representations on $R(\mathfrak{H})$ (which gives the analogue of the Weyl group) and on $A_\mathfrak{H}$.

**Theorem 6.18.** *Let $\mathfrak{H} \leq \mathrm{Aut}(G)$ have local tidy subgroups and let $\beta \in \mathrm{Aut}(G)$ satisfy $\beta \mathfrak{H} \beta^{-1} = \mathfrak{H}$. Then the automorphism $\alpha \mapsto \beta \alpha \beta^{-1}$ of $\mathfrak{H}$ induces a permutation $\phi \mapsto \beta.\phi$ of $\Phi$. Denote the normaliser of $\mathfrak{H}$ by $N_\mathfrak{H}$. The map*

$$\pi : N_\mathfrak{H} \ni \beta \mapsto (\phi \mapsto \beta.\phi), \ (N_\mathfrak{H} \to S_\Phi),$$

*is a homomorphism and $\mathfrak{H}$ is contained in its kernel.*

*The homomorphism $\pi$ induces a homomorphism $\tilde{\pi} : N_\mathfrak{H} \to \mathrm{Aut}\left(\bigoplus_{\phi \in \Phi} \mathbf{Z}\right)$ defined by $\tilde{\pi}(z_\phi) = (z_{\beta.\phi})$, $((z_\phi) \in \bigoplus_{\phi \in \Phi} \mathbf{Z})$ and $\tilde{\pi}(N_\mathfrak{H})$ leaves the rank group $R(\mathfrak{H})$ invariant. Hence $\tilde{\pi}$ induces representations of $N_\mathfrak{H}$ on the rank group, $R(\mathfrak{H})$, and on the corank group, $A_\mathfrak{H}$.*

**Proof.** Let $\phi$ be in $\Phi$. Then, for each finitely generated subgroup $\mathfrak{F}$ of $\mathfrak{H}$, there is a finitely generated $\mathfrak{L} \geq \mathfrak{F}$ and a $U$-eigenfactor, $V$, for $\mathfrak{L}$ such that

$$\phi(\alpha) = \begin{cases} [\alpha(V) : V], & \text{if } \alpha(V) \geq V \\ [V : \alpha(V)]^{-1}, & \text{if } \alpha(V) \leq V. \end{cases}$$

The homomorphism $\beta.\phi : \mathfrak{H} \to \mathbf{Q}^+$ defined by $\beta.\phi(\alpha) = \phi(\beta^{-1}\alpha\beta)$, $(\alpha \in \mathfrak{H})$, satisfies

$$\beta.\phi(\alpha) = \begin{cases} [\alpha(\beta(V)) : \beta(V)], & \text{if } \alpha(\beta(V)) \geq \beta(V) \\ [\beta(V) : \alpha(\beta(V))]^{-1}, & \text{if } \alpha(\beta(V)) \leq \beta(V), \end{cases}$$

where $\beta(V)$ is a $\beta(U)$-eigenfactor for $\beta\mathfrak{L}\beta^{-1} \geq \beta\mathfrak{F}\beta^{-1}$. Since $\beta\mathfrak{L}\beta^{-1}$ and $\beta\mathfrak{F}\beta^{-1}$ are finitely generated subgroups of $\mathfrak{H}$ it follows that $\beta.\phi$ is in $\Phi$ and the map $\pi(\beta) : \phi \mapsto \beta.\phi$ is a permutation of $\Phi$.

It is easily checked that $\pi$ is a homomorphism and that $\mathfrak{H}$ is contained in its kernel follows from Theorem 4.14. The remaining assertions are also easily verified. $\square$

*Remark* 6.19. In the course of the proof of Theorem 6.8 a set, $M_\mathfrak{H}$, of $n$-tuples of integers was defined, where $n$ is the rank of the finitely generated group $\mathfrak{H}$. This set is not an invariant for $\mathfrak{H}$ because it depends on a choice of isomorphism between $\mathfrak{H}/\mathfrak{N}$ and $\mathbf{Z}^n$. The number of points in $M_\mathfrak{H}$ is the factor number and is an invariant for $\mathfrak{H}$.

There are many geometrical invariants associated with $M_\mathfrak{H}$ as well. These include the number of extreme points in the convex hull of $M_\mathfrak{H}$ (as a subset of $\mathbf{R}^n$) and the volume of the convex hull. It may be that some of these invariants play a role in the structure theory of totally disconnected locally compact groups.

### Acknowledgement

The final form of this paper has been improved very much by discussions with Udo Baumgartner, Helge Glöckner and Jacqui Ramagge and I am grateful for their time and comments.

School of Mathematical and Physical Sciences, University of Newcastle, Callaghan, NSW 2308, Australia

george@frey.newcastle.edu.au